\documentclass[12pt]{article}

\usepackage{amsmath}
\usepackage{amssymb}
\usepackage{graphicx}
\usepackage{setspace}
\usepackage{amsfonts}
\usepackage{color}
\onehalfspace

\newtheorem{theorem}{Theorem}[section]

\newtheorem{definition}{Definition}[section]
\newtheorem{proposition}{Proposition}[section]

\begin{document}

\title{The Mahalanobis distance for functional data \\ with applications to classification}
\date{}
\author{Esdras Joseph, Pedro Galeano and Rosa E. Lillo\\
Departamento de Estad\'{\i}stica \\
Universidad Carlos III de Madrid}
\maketitle

\begin{abstract}
\noindent

This paper presents a general notion of Mahalanobis distance for functional data that extends the classical multivariate concept to situations where the observed data are points belonging to curves generated by a
stochastic process. More precisely, a new semi-distance for functional observations that generalize the usual Mahalanobis distance for multivariate datasets is introduced. For that, the development uses a regularized
square root inverse operator in Hilbert spaces. Some of the main characteristics of the functional Mahalanobis semi-distance are shown. Afterwards, new versions of several well known functional classification procedures
are developed using the Mahalanobis distance for functional data as a measure of proximity between functional observations. The performance of several well known functional classification procedures are compared with
those methods used in conjunction with the Mahalanobis distance for functional data, with positive results, through a Monte Carlo study and the analysis of two real data examples.


\end{abstract}

\newpage

\section{Introduction}

At the present time, there are a number of situations in different fields of applied sciences such as chemometrics, economics, image analysis, medicine, meteorology and speech recognition, among others, where it can be assumed
that the observed data are points belonging to functions defined over a given set. Functional data analysis (FDA) deals with such kind of observations. In practice, the values of the functions are available only at a finite
number of points and, as a general rule, functional samples may contain less functions than evaluation points. For this and other reasons, the majority of known multivariate tools can not be used for statistical analysis with this type of data since, by its nature, requires a different type of treatment. There are several methodologies for FDA being the most popular the one based on the use of basis functions such as Fourier and splines, see
Ramsay and Silverman (2005). Alternatively, other procedures, such as the nonparametric approach proposed by Ferraty and Vieu (2006), do not require the knowledge of the explicit form of the functions. The ideas developed in
this paper can be adapted to any of these situations. However, for easiness in exposition, the focus of this paper is on the basis functions approach.

Even if usual multivariate methods are not usually well suited for functional datasets, many multivariate techniques have inspired advances in FDA. The introduction of the notion of distance for functional data represents an
example. Usually, it is assumed that the set of functions has been generated by a functional random variable defined on a Hilbert space endowed with a certain distance. However, in the recent literature on functional data, there is little reference to the role played by distances between functional data, with the book of Ferraty and Vieu (2006) an exception. These authors have proposed semi-metrics well adapted for sample functions, including
semi-metrics based on functional principal components (FPC), partial least-squares (PLS) type semi-metrics and semi-metrics based on derivatives. However, common distances frequently used in multivariate data analysis
such as the Mahalanobis distance proposed by Mahalanobis (1936) have not been extended to the functional framework. The first contribution of this paper is to fill this gap and presents the funcional Mahalanobis semi-distance
that extends the multivariate Mahalanobis distance to the functional setting.

The use of distances in multivariate analysis is important in many different problems including classification, clustering, hypothesis testing and outlier detection, among others. In particular, several of the most well known
methods for classification analysis are distance-based. Under a functional perspective, the aim of classification procedures is to decide whether a function $\chi_{0}$ generated from a functional random variable $\chi $ belongs to one of $G$ classes using the information provided by $G$ independent training samples $\chi _{g1},\ldots ,\chi _{gn_{g}}$, where $g=1,\ldots ,G$. Here $\chi _{gi}$, for $i=1,\ldots ,n_{g}$, are independent replications of the
functional random variable $\chi $, measured on $n_{g}$ randomly chosen individuals from class $g$. Using this information, a functional classification method provides a classification rule that can be used to classify $\chi_{0}$. Nowadays, there is a wide variety of methods developed to solve this problem. For instance, several papers have proposed to classify functional observations by means of the functional principal
component scores. For instance, Hall et al. (2001) proposed a method that consists in obtain the functional principal component scores of the training samples, then estimate nonparametrically the probability densities of the sets of functional principal component scores and finally estimate the posterior probability that $\chi_{0}$ is of a given class using the Bayes classification rule. This approach was considered by Glendinning and
Herbert (2003) for shape classification. Under a similar perspective, Leng and M\"{u}ller (2006) proposed a method of classifying collections of temporal gene expression curves by means of functional logistic regression on the
functional principal component scores of the training samples. Also, Song et al. (2008) compared several multivariate classification methods on the the basis expansion coefficients of the training samples for classifying
time-course gene expression data. On the other hand, the popular nearest neighbor classification rule has been also extended to functional data. For instance, Biau et al. (2005) proposed to filter the training samples in the Fourier basis and to apply the kNN method to the first Fourier coefficients of the expansion, while Ba\'illo et al. (2011) derived several consistency results of the kNN procedure for a particular type of Gaussian processes.
Additionally, the centroid method based on assign the function to the group with closer mean has been adapted to the functional framework by Delaigle and Hall (2012). Alternatively, several papers have extended the Fisher's
discriminant analysis to the functional framework. The idea of these methods is to project the observations into a finite dimensional space where the classes are separated as much as possible. The transformed functions are called discriminant functions. Then, the new function $\chi_{0}$ is also projected in this space and it is classified using the Bayes classification rule. In particular, James and Hastie (2001) used a natural cubic spline basis plus
random error to model the observations from each individual. The spline is parameterized using a basis function multiplied by a coefficient vector, that is modeled using a Gaussian distribution. The observed functions can then be pooled to estimate the mean and covariance for each class by means of an Expectation-Maximization (EM) algorithm that are used to obtain the discriminant functions. Alternatively, Preda et al. (2007) used functional PLS regression to obtain the discriminant functions while Shin (2008) considered an approach based on reproducing kernel Hilbert spaces. Finally, Ferraty and Vieu (2003) have proposed a method based on estimating nonparametrically the posterior probability that the new function $\chi_{0}$ is of a given class, L\'{o}pez-Pintado and Romo (2006), Cuevas et al. (2007) and Sguera et al. (2012) have proposed classifiers based on the notion of data depth that are well suited for datasets containing outliers, Rossi and Villa (2006) and Martin-Barragan et al. (2013) have investigated the use of support vector machines (SVMs) for functional data, Wang et al. (2007) have considered classification for functional data by Bayesian modeling with wavelet basis functions, Epifanio (2008) has developed classifiers based on shape descriptors, Araki et al. (2009) have considered functional logistic classification, and, finally, Alonso et al. (2012) have proposed a weighted distance approach. Note that, when a distance is required, these papers use the $L^1$, $L^2$ and $L^\infty$ distances which are well defined in Hilbert spaces. The second contribution of this paper is to show that several simple classification procedures including the kNN procedure, the centroid method and functional Bayes classification rules can be used in conjunction with the functional Mahalanobis semi-distance as the criterion of proximity between functions to get very good classification rates without the need of much higher sophisticated classification methods. Several Monte Carlo experiments suggest that methods based on the functional Mahalanobis semi-distance leads to better classification rates than other alternatives.

The rest of this paper is organized as follows. Section $2$ introduces the functional Mahalanobis semi-distance and shows some of its main characteristics. Section $3$ reviews several classification methods for functional data and provides new approaches to these methods based on the functional Mahalanobis semi-distance. Section $4$ analyzes the empirical properties of the procedures via several Monte Carlo experiments and illustrates the good
behavior of the classification methods in conjunction with the functional Mahalanobis semi-distance through of the analysis of two real data examples. Finally, some conclusions are drawn in Section $5$.

\section{The funcional Mahalanobis semi-distance}

\subsection{Definitions and some characteristics}

This section presents the functional Mahalanobis semi-distance that generalizes the Mahalanobis distance for multivariate random variables to the functional framework. Let $\mathbf{x}$ be a multivariate continuous random
variable defined in $\mathbb{R}^{p}$ with mean vector $\mathbf{m}_{\mathbf{x}}=E\left[ \mathbf{x}\right] $ and definite positive covariance matrix $\mathbf{C}_{\mathbf{x}}=E\left[ \left( \mathbf{x}-\mathbf{m}_{\mathbf{x}}\right) \left( \mathbf{x}-\mathbf{m}_{\mathbf{x}}\right) ^{\prime }\right] $. The Mahalanobis distance between the random variable $\mathbf{x}$ and its mean vector $\mathbf{m}_{\mathbf{x}}$ is the Euclidean norm of the random vector $\mathbf{C}_{\mathbf{x}}^{-1/2}\left( \mathbf{x}-\mathbf{m}_{\mathbf{x}}\right) $ that can be written (see, Mahalanobis, 1936) as:
\begin{equation}
\begin{tabular}{c}
$d_{M}\left( \mathbf{x},\mathbf{m}_{\mathbf{x}}\right) =\left\Vert \mathbf{C}_{\mathbf{x}}^{-1/2}\left( \mathbf{x}-\mathbf{m}_{\mathbf{x}}\right)\right\Vert _{E}=$ \\
$=\left\langle \mathbf{C}_{\mathbf{x}}^{-1/2}\left( \mathbf{x}-\mathbf{m}_{\mathbf{x}}\right) ,\mathbf{C}_{\mathbf{x}}^{-1/2}\left( \mathbf{x}-\mathbf{m}_{\mathbf{x}}\right) \right\rangle _{E}^{1/2}=\left[ \left( \mathbf{x}-\mathbf{m}_{\mathbf{x}}\right) ^{\prime }\mathbf{C}_{\mathbf{x}}^{-1}\left(\mathbf{x}-\mathbf{m}_{\mathbf{x}}\right) \right] ^{1/2},$
\end{tabular}
\label{Pred}
\end{equation}
where $\left\Vert \mathbf{\cdot }\right\Vert _{E}$ and $\left\langle \mathbf{\cdot },\mathbf{\cdot }\right\rangle _{E}$ denote the Euclidean norm and the usual inner product in $\mathbb{R}^{p}$, respectively. The main characteristic of the multivariate Mahalanobis distance is that it takes into account the correlation structure of the multivariate random variable $\mathbf{x}$. Moreover, the multivariate Mahalanobis distance is scale
invariant. For future developments, it is important to note that the Mahalanobis distance can be written in terms of the principal component scores of $\mathbf{x}$. For that, let $\mathbf{v}_{1},\ldots ,\mathbf{v}_{p}$
be the eigenvectors of the covariance matrix $\mathbf{C}_{\mathbf{x}}$ associated with positive eigenvalues $a_{1}\geq \cdots \geq a_{p}>0$, and let $\mathbf{V}$ be the $p\times p$ matrix whose columns are the
eigenvectors of the covariance matrix $\mathbf{C}_{\mathbf{x}}$, i.e., $\mathbf{V=}\left[ \mathbf{v}_{1}|\cdots |\mathbf{v}_{p}\right] $. Then, the vector of principal component scores given by $\mathbf{s}=\mathbf{V}^{\prime }\left( \mathbf{x}-\mathbf{m}_{\mathbf{x}}\right) $, is a multivariate random variable with zero mean vector and diagonal covariance matrix. As a consequence, $\mathbf{x}$ can be written in terms of the principal component scores in the following way:
\begin{equation}
\mathbf{x}=\mathbf{m}_{\mathbf{x}}+\mathbf{Vs.}  \label{xmu}
\end{equation}%
On the other hand, the singular value decomposition of $\mathbf{C}_{\mathbf{x}}$, i.e., $\mathbf{C}_{\mathbf{x}}=\mathbf{V}\mathbf{A}\mathbf{V}^{\prime }$, where $\mathbf{A}$ is a diagonal matrix with
the ordered eigenvalues $a_{1},\ldots ,a_{p}$ in the main diagonal, allows to write the inverse of $\mathbf{C}_{\mathbf{x}}$ in terms of $\mathbf{V}$ and $\mathbf{A}$ as follows:
\begin{equation}
\mathbf{C}_{\mathbf{x}}^{-1}=\mathbf{VA }^{-1}\mathbf{V}^{\prime}.  \label{sinv}
\end{equation}
Now, (\ref{xmu}) and (\ref{sinv}) leads to the following expression of the Mahalanobis distance between the random variable $\mathbf{x}$ and its mean vector $\mathbf{m}_{\mathbf{x}}$ in terms of the principal component
scores:
\begin{equation}
d_{M}\left( \mathbf{x},\mathbf{m}_{\mathbf{x}}\right) =\left( \mathbf{%
s}^{\prime }\mathbf{V}^{\prime }\mathbf{VA }^{-1}\mathbf{V}^{\prime }%
\mathbf{Vs}\right) ^{1/2}=\left( \mathbf{s}^{\prime }\mathbf{A }^{-1}%
\mathbf{s}\right) ^{1/2}=\left( \mathbf{z}^{\prime }\mathbf{z}\right) ^{1/2},
\label{scorescomp}
\end{equation}
where $\mathbf{z=A}^{-1/2}\mathbf{s}$ is the random vector of standardized principal component scores. In other words, the Mahalanobis distance between $\mathbf{x}$ and $\mathbf{m}_{\mathbf{x}}$ can be written as the Euclidean norm of the standardized principal component scores.

As mentioned before, the main goal of this section is to generalize the multivariate Mahalanobis distance to the functional setting. However, the proposal does not lead to a functional distance but to a functional semi-distance. The reasons of this will be clear once the functional Mahalanobis semi-distance is presented. For that, let $\chi $ be a functional random variable defined in the infinite dimensional space $L^{2}(T)$, i.e., the space of squared integrable functions in the closed interval $T=\left[ a,b\right] $ of the real line. It is assumed that the functional random variable $\chi $ has a functional mean $\mu _{\chi }(t)=E[\chi (t)]$ and a covariance operator $\Gamma _{\chi }$ given by:
\begin{equation}
\Gamma _{\chi }(\eta )=E[(\chi -\mu _{\chi })\otimes (\chi -\mu _{\chi
})(\eta )],  \label{cova}
\end{equation}%
such that, for any $\eta \in L^{2}(T)$,
\begin{equation}
(\chi -\mu _{\chi })\otimes (\chi -\mu _{\chi })(\eta )=\left\langle \chi
-\mu _{\chi },\eta \right\rangle (\chi -\mu _{\chi }),  \label{prod1}
\end{equation}%
where $\left\langle .,.\right\rangle $ denotes the usual inner product on $L^{2}(T)$, i.e.:
\[
\left\langle \chi -\mu _{\chi },\eta \right\rangle =\int_{T}\left( \chi
(t)-\mu _{\chi }(t)\right) \eta (t)dt.
\]%
The covariance operator $\Gamma _{\chi }$ in (\ref{cova}) is a well-defined compact operator so long as $E\left[ \left\Vert \chi \right\Vert _{2}^{4}\right] <\infty $ (see Hall and Hosseini-Nasab, 2006), where $\left\Vert .\right\Vert _{2}$ denotes the usual norm in $L^{2}(T)$. Under this assumption, there exists a sequence of non-negative eigenvalues of $\Gamma _{\chi }$, denoted by $\lambda _{1}\geq \lambda _{2}\geq \cdots$, where $\sum_{k=1}^{\infty }\lambda_{k}<\infty $, and a set of orthonormal eigenfunctions of $\Gamma _{\chi }$, denoted by $\psi _{1},\psi _{2},\ldots $ such that $\Gamma _{\chi }(\psi_{k})=\lambda _{k}\psi _{k}$, for $k=1,2,\ldots $ The eigenfunctions $\psi_{1},\psi _{2},\ldots $ form an orthonormal basis in $L^{2}(T)$ and allows to write the Karhunen-Lo\`{e}ve expansion of the functional random variable $\chi $ (see Hall and Housseini-Nassab (2006)), in terms of the elements of the basis as follows:
\begin{equation}
\chi =\mu _{\chi }+\sum_{k=1}^{\infty }\theta _{k}\psi _{k},  \label{loeve23}
\end{equation}
where $\theta _{k}=\left\langle \chi -\mu _{\chi },\psi _{k}\right\rangle $, for $k=1,2,\ldots $ are the functional principal component scores of $\chi $. It is well known that the functional principal component scores $\theta _{k}$, for $k=1,2,\ldots$ are uncorrelated random variables with zero mean and variance $\lambda_{k}$ since $\psi _{1},\psi _{2},\ldots $ are orthonormal.

In order to obtain a similar expression to (\ref{Pred}) in the functional setting, it is necessary to define the inverse of the covariance operator, $\Gamma _{\chi }^{-1}$. It exists under certain circumstances. However, even in this case, $\Gamma _{\chi }^{-1}$ is unbounded and not continuous. Mas (2007) has proposed a regularized inverse operator which is a linear operator \textquotedblleft close\textquotedblright\ to $\Gamma _{\chi }^{-1}$ and having good properties. For that, if $\Gamma _{\chi }^{-1}$ exists, this is given by:
\[
\Gamma _{\chi }^{-1}(\zeta )=\sum_{k=1}^{\infty }\frac{1}{\lambda _{k}}(\psi _{k}\otimes \psi _{k})(\zeta ),
\]%
where $\zeta $ is a function in the range of $\Gamma _{\chi }$. Then, the regularized inverse operator, denoted by $\Gamma _{K}^{-1}$, is defined as:
\[
\Gamma _{K}^{-1}(\zeta )=\sum_{k=1}^{K}\frac{1}{\lambda _{k}}(\psi_{k}\otimes \psi _{k})(\zeta ),
\]%
where $K$ is a given threshold. Similarly, it is also possible to give a regularized square root inverse operator given by:
\begin{equation}
\Gamma _{K}^{-1/2}(\zeta )=\sum_{k=1}^{K}\frac{1}{\lambda _{k}^{1/2}}(\psi_{k}\otimes \psi _{k})(\zeta ),  \label{pro12}
\end{equation}%
that allows to define the functional Mahalanobis semi-distance between $\chi$ and $\mu _{\chi }$ inspired on (\ref{Pred}) as follows:

\begin{definition}
Let $\chi $ be a functional random variable defined in $L^{2}(T)$ with mean function $\mu _{\chi }$ and compact covariance operator $\Gamma _{\chi }$. The Mahalanobis semi-distance between $\chi $ and $\mu _{\chi }$, denoted by
$d_{FM}^{K}(\chi ,\mu _{\chi })$, is defined as:
\[
d_{FM}^{K}(\chi ,\mu _{\chi })=\left\langle \Gamma _{K}^{-1/2}(\chi -\mu_{\chi }),\Gamma _{K}^{-1/2}(\chi -\mu _{\chi })\right\rangle ^{1/2}.
\]
\end{definition}

As noted before, the multivariate Mahalanobis distance may be expressed in terms of the principal component scores of the multivariate random variable $\mathbf{x}$. Similarly, it is possible to express the functional Mahalanobis
semi-distance in terms of the functional principal component scores of the functional random variable $\chi $ as stated in the next proposition, that is proved in the appendix:

\begin{proposition}
The functional Mahalanobis semi-distance between $\chi $ and $\mu _{\chi }$ can be written as follows:
\begin{equation}
d_{FM}^{K}(\chi ,\mu _{\chi })=\left( \sum_{k=1}^{K}\omega_{k}^{2}\right) ^{1/2},  \label{dmk}
\end{equation}
where $\omega _{k}=\theta _{k}/\lambda _{k}^{1/2}$, for $k=1,\ldots ,K$, are the standardized functional principal component scores.
\end{proposition}

Therefore, as in the multivariate case, the functional Mahalanobis semi-distance between $\chi$ and $\mu _{\chi }$ is the Euclidean norm of the standardized functional principal component scores. This property provides a simple
way to compute the functional Mahalanobis semi-distance in practice. It is also interesting to extend the definition of functional Mahalanobis semi-distance to the general situation of distance between two independent and
identically distributed functional random variables.

\begin{definition}
Let $\chi _{1}$ and $\chi _{2}$ be two functional random variables defined in $L^{2}(T)$ independent and identically distributed with mean function $\mu _{\chi }$ and compact covariance operator $\Gamma _{\chi }$. The
functional Mahalanobis semi-distance between the functions $\chi _{1}$ and $\chi _{2}$, denoted by $d_{FM}^{K}(\chi _{1},\chi _{2})$, is given by:
\[
d_{FM}^{K}(\chi _{1},\chi _{2})=\left\langle \Gamma _{K}^{-1/2}(\chi_{1}-\chi _{2}),\Gamma _{K}^{-1/2}(\chi _{1}-\chi _{2})\right\rangle ^{1/2}.
\]
\end{definition}

The previous definition leads to the following proposition proven in the appendix:

\begin{proposition}
The functional Mahalanobis semi-distance between $\chi _{1}$ and $\chi _{2}$ can be written as follow:
\begin{equation}
d_{FM}^{K}(\chi _{1},\chi _{2})=\left( \sum_{k=1}^{K}\left( \omega_{1k}-\omega _{2k}\right) ^{2}\right) ^{1/2},  \label{mah1}
\end{equation}%
where $\omega _{1k}=\theta _{1k}/\lambda _{k}^{1/2}$ and $\omega_{2k}=\theta _{2k}/\lambda _{k}^{1/2}$, for $k=1,2,\ldots $ are the standardized functional principal component scores of $\chi _{1}$ and $\chi_{2}$, respectively.
\end{proposition}

Therefore, the functional Mahalanobis semi-distance between two independent and identically distributed functional random variables can be written as the Euclidean distance between the standardized functional principal component scores of both functional random variables. The next result shows that $d_{FM}^{K}$ is indeed a functional semi-distance.

\begin{proposition}
Let $\chi _{1}$, $\chi _{2}$ and $\chi _{3}$ be three independent and identically distributed functional random variables defined in $L^{2}(T)$ with mean function $\mu _{\chi }$ and compact covariance operator $\Gamma _{\chi
} $. For any positive integer $K$, $d_{FM}^{K}$ verifies the following three properties:

\begin{enumerate}
\item $d_{FM}^{K}(\chi _{1},\chi _{2})\geq 0$.

\item $d_{FM}^{K}(\chi _{1},\chi _{2})=d_{FM}^{K}(\chi _{2},\chi _{1})$.

\item $d_{FM}^{K}(\chi _{1},\chi _{2})\leq d_{FM}^{K}(\chi _{1},\chi
_{3})+d_{FM}^{K}(\chi _{3},\chi _{2})$.
\end{enumerate}

Consequently, $d_{FM}^{K}$ is a functional semi-distance.
\end{proposition}

It is well known that if the multivariate random variable $\mathbf{x}$ has a $p$-dimensional Gaussian distribution, then it is easy to see that $d_{M}^{2}(\mathbf{x},\mathbf{m}_{\mathbf{x}})$ has a $\chi _{p}^{2}$
distribution and, consequently, $E\left[ d_{M}^{2}(\mathbf{x},\mathbf{m}_{\mathbf{x}})\right] =p$ and $V\left[ d_{M}^{2}(\mathbf{x},\mathbf{m}_{\mathbf{x}})\right] =2p$. To end this section, the following theorem shows
a similar result for the functional Mahalanobis semi-distance.

\begin{theorem}
If $\chi $ is a Gaussian process, $d_{FM}^{K}(\chi ,\mu _{\chi })^{2}\sim\chi _{K}^{2}$, so that $E\left[ d_{FM}^{K}(\chi ,\mu _{\chi })^{2}\right]=K $ and $V\left[ d_{FM}^{K}(\chi ,\mu _{\chi })^{2}\right] =2K$.
\end{theorem}

\subsection{Practical implementation}

In practice, the functions are not observed continuously over all the points in the closed interval $T=\left[ a,b\right] $, so that calculation of the functional Mahalanobis semi-distances as defined in (\ref{dmk}) and (\ref{mah1}) is not possible. Assume now that a dataset is observed with the following form:
\begin{equation}
\left\{ \chi _{i}\left( t_{i,j}\right) :i=1,\ldots ,n\text{ and }j=1,\ldots
,J_{i}\right\} ,  \label{functionaldataset}
\end{equation}
where $n$ is the number of observed curves and $J_{i}$ is the number of observations of the function $\chi _{i}$ at the points $t_{i,1},\ldots,t_{i,J_{i}}$. Note that it is not assumed that the observation points are
the same for all the functions not even their numbers. In this situation, the usual approach to obtain closed form expressions of the set of functions is to use basis functions. In general, a basis is a system of functions,
denoted by $\phi _{m}$, for $m=1,2,\ldots $, orthogonal or not, such that, for $i=1,\ldots ,n$:
\[
\chi _{i}\left( t\right) \simeq \sum\limits_{m=1}^{M}\beta_{im}\phi _{m}\left(t\right) ,
\]
where $\beta_{im}$, for $m=1,\ldots ,M$, are the coefficients of the expansion. The number of basis functions, $M$, should be chosen on a case by case basis, although, $M$ is usually chosen such that the functional approximations are close to the original counterparts with some smoothing that eliminates the most obvious noise. The choice of the basis is also important. There are several possibilities including polynomial, wavelets, Fourier and splines basis, among others. For periodic or nearly periodic datasets, Fourier basis is an adequate choice. For nonperiodic datasets, B-splines are typically used. See Ramsay and Silverman (2005) for more information on basis functions. The simplest method to effectively estimate the coefficients of the expansion is carried out by minimizing:
\[
\left( \sum_{j=1}^{J_{i}}\left[ \chi _{i}\left( t_{i,j}\right)
-\sum\limits_{m=1}^{M}\beta_{im}\phi _{m}\left( t_{i,j}\right) \right]
^{2}\right) ^{1/2}.
\]
Now, with the smoothed functional sample, it is possible to estimate the functional mean $\mu _{\chi }$ with the sample functional mean, $\widehat{\mu }_{\chi }$, given by:
\[
\widehat{\mu }_{\chi }=\frac{1}{n}\sum\limits_{i=1}^{n}\chi _{i},
\]
and the covariance operator $\Gamma _{\chi }$ with the sample covariance operator, $\widehat{\Gamma }_{\chi }\left( \eta \right)$, such that, for any $\eta \in L^{2}(T)$:
\[
\widehat{\Gamma }_{\chi }\left( \eta \right) =\frac{1}{n}\sum\limits_{i=1}^{n}\left\langle \chi_i -\widehat{\mu} _{\chi },\eta \right\rangle (\chi_i-\widehat{\mu} _{\chi }).
\]
Then, eigenfunctions and eigenvalues of the covariance operator $\Gamma_{\chi }$ can be approximated with those of $\widehat{\Gamma }_{\chi }$ leading to estimates $\widehat{\psi }_{1},\widehat{\psi }_{2},\ldots $ and $\widehat{\lambda }_{1},\widehat{\lambda }_{2},\ldots$ respectively. Therefore, the functional principal component scores corresponding to curve $\chi _{i}$, i.e., $\theta _{i,k}=\left\langle \chi _{i}-\mu _{\chi },\psi
_{k}\right\rangle $, are estimated with $\widehat{\theta }_{i,k}=\left\langle \chi _{i}-\widehat{\mu }_{\chi },\widehat{\psi }_{k}\right\rangle$, for $k=1,2,\ldots $ that allows us to define the functional Mahalanobis semi-distance between $\chi _{i}$ and the functional sample mean $\widehat{\mu }_{\chi }$ as follows:
\[
d_{FM}^{K}(\chi _{i},\widehat{\mu }_{\chi })=\left( \sum_{k=1}^{K}\widehat{\omega }_{ik}^{2}\right) ^{1/2},
\]
where $\widehat{\omega }_{ik}=\widehat{\theta }_{i,k}/\widehat{\lambda }_{k}^{1/2}$, for $k=1,\ldots ,K$, are the sample standardized functional principal component scores. Similarly, the functional Mahalanobis semi-distance
between two functions of the sample, $\chi _{i}$ and $\chi _{i^{\prime }}$, can be written as follows:
\[
d_{FM}^{K}(\chi _{i},\chi _{i^{\prime }})=\left( \sum_{k=1}^{K}\left(
\widehat{\omega }_{ik}-\widehat{\omega }_{i^{\prime }k}\right) ^{2}\right)
^{1/2},
\]
where $\widehat{\omega }_{i^{\prime }k}=\widehat{\theta }_{i^{\prime },k}/\widehat{\lambda }_{k}^{1/2}$, for $k=1,\ldots ,K$.

\section{Classification with the functional Mahalanobis semi-distance}

Among all the possible applications of the functional Mahalanobis semi-distance introduced in the previous Section, this paper focuses in the supervised classification problem in the functional setting. Consider a sample
of functional observations such that it is known in advance that each function comes from one of $G$ predefined classes. Therefore, the whole sample can be split in $G$ subsamples, denoted by $\chi _{g1},\ldots ,\chi
_{gn_{g}}$, for $g=1,\ldots ,G$, respectively, where $n=n_{1}+\cdots +n_{G}$ is the sample size of the whole dataset. Then, the idea is to use the information provided by the set of observations to construct
classification rules that can be used to classify a new ungrouped functional observation $\chi _{0}$. The aim of this section is to propose new procedures based on the combination of well known functional classification
methods with the functional Mahalanobis semi-distance as a measure of proximity between functional objects. In particular, four procedures are presented.

\subsection{The k-nearest neighbor (kNN) procedure}

The k-nearest neighbor (kNN) procedure is one of the most popular methods used to perform supervised classification in multivariate settings. The method is very simple and appears to have a very good performance in many
situations. Its generalization to infinite-dimensional spaces has been studied by Biau et al. (2005), C\'{e}rou and Guyader (2006) and Ba\'{\i}llo et al. (2011), among others. The kNN method starts by computing the distances
between the new function to classify, $\chi _{0}$, and all the functions in the observed sample. Next, the method finds the $k$ functional observations in the sample closest in distance to $\chi _{0}$. Finally, the
new observation $\chi _{0}$ is classified using majority of votes among the $k$ neighbors. C\'{e}rou and Guyader (2006) have shown that the kNN procedure is not universally consistent. However, these authors have obtained
sufficient conditions for consistency of the kNN classifier when the functional random variable takes values in a separable metric space. Additionally, Ba\'{\i}llo et al. (2011) have shown that the optimal classification rule can be explicitly obtained for a class of Gaussian processes with triangular covariance operators. The previous papers have considered three functional distances for the kNN classifier: the $L^{1},$ $L^{2}$ and $L^{\infty }$
distances. In particular, the $L^{1},$ $L^{2}$ and $L^{\infty }$ distances between $\chi _{0}$ and the functional observation $\chi _{gi}$ for $g=1,\ldots ,G$ and $i=1,\ldots ,n_{g}$ are given by:
\begin{gather*}
d_{1}\left( \chi _{0},\chi _{gi}\right) =\int_{T}\left\vert \chi _{0}\left(
t\right) -\chi _{gi}\left( t\right) \right\vert dt, \\
d_{2}\left( \chi _{0},\chi _{gi}\right) =\left( \int_{T}\left( \chi
_{0}\left( t\right) -\chi _{gi}\left( t\right) \right) ^{2}dt\right) ^{1/2},
\end{gather*}
and,
\begin{equation*}
d_{\infty }\left( \chi _{0},\chi _{gi}\right) =\sup \left\{ \left\vert \chi
_{0}\left( t\right) -\chi _{gi}\left( t\right) \right\vert :t\in T\right\} ,
\end{equation*}
respectively. Note that in order to compute the $L^{1},$ $L^{2}$ and $L^{\infty }$ distances it is necessary to first smooth the discretized values of the function $\chi _{0}$ as seen in Section 2.2. Also, it is important to note that no information about the class membership is used to compute the previous distances.

On the other hand, the kNN classifier can be used in conjunction with the functional Mahalanobis semi-distance. Contrary to the previous distances, two different ways to compute the functional Mahalanobis semi-distance in classification problems are in order. In a first case, assume that the functional means under class $g$, denoted by $\mu _{\chi_{g}}$, are different but the covariance operator, denoted by $\Gamma _{\chi}$, is the same for all the classes. Then, the functional means, $\mu _{\chi_{g}}$, are estimated using the functional sample mean of the functions in class $g$, i.e.:
\begin{equation}
\widehat{\mu }_{\chi _{g}}=\frac{1}{n_{g}}\sum\limits_{i=1}^{n_{g}}\chi_{gi},  \label{MeanEst}
\end{equation}
while the common covariance operator, $\Gamma _{\chi }$, is estimated with the within class covariance operator given by:
\begin{equation}
\widehat{\Gamma }_{\chi }\left( \eta \right) =\frac{1}{n}\sum\limits_{g=1}^{G}\sum\limits_{i=1}^{n_{g}}\left\langle \chi _{gi}-\widehat{\mu }_{\chi _{g}},\eta \right\rangle \left( \chi _{gi}-\widehat{\mu }_{\chi
_{g}}\right) ,  \label{CovEst}
\end{equation}%
for $\eta \in L^{2}\left( T\right) $. Now, the functional Mahalanobis semi-distance between $\chi _{0}$ and the functional observation $\chi _{gi}$ for $g=1,\ldots ,G$ and $i=1,\ldots ,n_{g}$ is given by:
\begin{equation}
d_{FM}^{K}\left( \chi _{0},\chi _{gi}\right) =\left(
\sum\limits_{k=1}^{K}\left( \widehat{\omega }_{g0k}-\widehat{\omega }%
_{gik}\right) ^{2}\right) ^{1/2},  \label{dFMSim}
\end{equation}%
where $\widehat{\omega }_{g0k}=\widehat{\theta }_{g0k}/\widehat{\lambda }_{k}^{1/2}$ and $\widehat{\omega }_{gik}=\widehat{\theta }_{gik}/\widehat{\lambda }_{k}^{1/2}$, respectively, are the standardized sample functional
principal component scores given by $\widehat{\theta }_{g0k}=\left\langle \chi _{0}-\widehat{\mu }_{\chi _{g}},\widehat{\psi }_{k}\right\rangle $ and $\widehat{\theta }_{gik}=\left\langle \chi _{gi}-\widehat{\mu }_{\chi _{g}},\widehat{\psi }_{k}\right\rangle $, respectively. Here, $\widehat{\psi }_{1},\ldots ,\widehat{\psi }_{K}$ and $\widehat{\lambda }_{1},\ldots ,\widehat{\lambda }_{K}$ are the eigenfunctions and eigenvalues of the sample
within class covariance operator (\ref{CovEst}). Similarly, the functional principal components (FPC) semi-distance proposed by Ferraty and Vieu (2006) between $\chi _{0}$ and the functional observation $\chi _{gi}$ for
$g=1,\ldots ,G$ and $i=1,\ldots ,n_{g}$, can be written as follows:
\begin{equation}
d_{FPC}^{K^{\prime}}\left( \chi _{0},\chi _{gi}\right) =\left(\sum\limits_{k=1}^{K^{\prime}}\left( \widehat{\theta }_{g0k}-\widehat{\theta }_{gik}\right) ^{2}\right) ^{1/2},  \label{FPC1}
\end{equation}
where $K^{\prime}$ is a certain threshold. In a second case, assume that both, the functional means and the covariance operators, denoted by $\Gamma _{\chi _{g}}$, are different for the classes $1,\ldots ,G$. Then, the functional means, $\mu _{\chi _{g}}$, are estimated using (\ref{MeanEst}), while the covariance operator of each class is estimated using the functional sample covariance operator of the functions in class $g$, i.e.:
\begin{equation}
\widehat{\Gamma }_{\chi _{g}}\left( \eta \right) =\frac{1}{n_{g}}\sum\limits_{i=1}^{n_{g}}\left\langle \chi _{gi}-\widehat{\mu }_{\chi_{g}},\eta \right\rangle \left( \chi _{gi}-\widehat{\mu }_{\chi
_{g}}\right) ,  \label{CovEst2}
\end{equation}
for $\eta \in L^{2}\left( T\right) $. Now, the functional Mahalanobis semi-distance between $\chi _{0}$ and the functional observation $\chi _{gi}$ for $g=1,\ldots ,G$ and $i=1,\ldots ,n_{g}$, is like in (\ref{dFMSim}) but
here $\widehat{\omega }_{g0k}=\widehat{\theta }_{g0k}/\widehat{\lambda }_{gk}^{1/2}$ and $\widehat{\omega }_{gik}=\widehat{\theta }_{gik}/\widehat{\lambda }_{gk}^{1/2}$, respectively, where $\widehat{\theta }_{g0k}=\left\langle
\chi _{0}-\widehat{\mu }_{\chi _{g}},\widehat{\psi }_{gk}\right\rangle $, $\widehat{\theta }_{gik}=\left\langle \chi _{gi}-\widehat{\mu }_{\chi _{g}},\widehat{\psi }_{gk}\right\rangle $ and $\widehat{\psi }_{g1},\ldots
,\widehat{\psi }_{gK}$ and $\widehat{\lambda }_{g1},\ldots ,\widehat{\lambda }_{gK}$ are the eigenfunctions and eigenvalues of the sample covariance operator in (\ref{CovEst2}), respectively. Also, the FPC semi-distance in this
second case can be written as in (\ref{FPC1}) but considering the same sample functional scores obtained with the eigenfunctions from the covariance operator (\ref{CovEst2}), as before.

\subsection{The centroid procedure}

The centroid procedure for functional datasets, proposed by Delaigle and Hall (2012), is probably the fastest and simplest classification method for functional observations. The centroid method consists in assigning a new function $\chi_{0}$ to the class with closer mean. Note that any functional distance can be used to implement the procedure. In particular, Delaigle and Hall (2012) considered the case of $G=2$ classes that have different mean and a common covariance operator and proposed to project the functions into a given direction and then compute the squared Euclidean distance between the observations. More precisely, Delaigle and Hall (2012) proposed to use the centroid classifier with the distance between $\chi _{0}$ and the sample functional mean $\widehat{\mu }_{\chi _{g}}$, for $g=1,2$, denoted by $DH$, and given by:
\begin{equation}
d_{DH}\left( \chi _{0},\widehat{\mu }_{\chi _{g}}\right) =\left\vert\sum\limits_{k=1}^{K^{\prime\prime}}\widehat{\omega }_{0gk}\widehat{\delta }
_{12k}\right\vert ,\label{DH}
\end{equation}
where $K^{\prime\prime}$ is a certain threshold, $\widehat{\omega }_{0gk}$ is computed as in the previous subsection assuming a common covariance operator and
\begin{equation*}
\widehat{\delta }_{12k}=\frac{\left\langle \widehat{\mu }_{\chi _{2}}-%
\widehat{\mu }_{\chi _{1}},\widehat{\psi }_{k}\right\rangle }{\widehat{%
\lambda }_{k}^{1/2}},
\end{equation*}
for $k=1,\ldots ,K^{\prime\prime}$.

Of course, other distances can be applied in the general case of $G$ classes. In particular, the $L^{1},$ $L^{2}$ and $L^{\infty }$ distances and the two versions of the functional Mahalanobis and functional principal components semi-distances introduced before, between $\chi _{0}$ and the sample functional mean $\widehat{\mu }_{\chi_{g}}$, for $g=1,\ldots ,G$, can be used. In particular, the semi-distances are computed similarly in the
previous Section but replacing $\chi _{gi}$ with $\widehat{\mu }_{\chi _{g}}$.

\subsection{The functional linear and quadratic Bayes classification rules}

In multivariate statistics, the Bayes classification rule is derived as follows. Let $\boldsymbol{x}$ be a $p$-dimensional continuous random variable and let $f_{1},\ldots ,f_{G}$ be the corresponding density functions of
$\boldsymbol{x}$ under the $G$ classes. Let $\pi _{1},\ldots ,\pi _{G}$ be the prior probabilities assigned to the $G$ classes, verifying $\pi _{1}+\cdots +\pi_{G}=1$. Using the Bayes Theorem, the posterior probability that a
new observation $\boldsymbol{x_{0}}$ generated from $\boldsymbol{x}$ comes from class $g$ is given by:
\begin{equation}
P(g|\boldsymbol{x_{0}})=\frac{\pi _{g}f_{g}(\boldsymbol{x_{0}})}{\pi
_{1}f_{1}(\boldsymbol{x_{0}})+\cdots +\pi _{G}f_{G}(\boldsymbol{x_{0}})}, \label{PP}
\end{equation}
respectively, where $P(1|\boldsymbol{x_{0}})+\cdots +P(G|\boldsymbol{x_{0}})=1$. The Bayes rule classifies $\boldsymbol{x_{0}}$ in the class with largest posterior probability. In other words, $\boldsymbol{x_{0}}$ is
classified in class $g$ if $\pi _{g}f_{g}(\boldsymbol{x_{0}})$ is maximum. In particular, if the $f_{g}$ densities are assumed to be Gaussian with different means $\mathbf{m}_{\boldsymbol{x}_{g}}$ but identical covariance
matrix $\mathbf{C}_{\boldsymbol{x}}$, this is equivalent to classify $\boldsymbol{x_{0}}$ in class $g$ if:
\begin{equation*}
d_{M}\left( \boldsymbol{x_{0}},\mathbf{m}_{\boldsymbol{x}_{g}}\right)^{2}
-2\log \pi _{g}
\end{equation*}
is minimum, where $d_{M}\left( \boldsymbol{x_{0}},\mathbf{m}_{\boldsymbol{x}_{g}}\right)^{2} =\left( \boldsymbol{x_{0}}-\mathbf{m}_{\boldsymbol{x}_{g}}\right) ^{\prime }\mathbf{C}_{\boldsymbol{x}}^{-1}\left(
\boldsymbol{x_{0}}-\mathbf{m}_{\boldsymbol{x}_{g}}\right) $ is the squared Mahalanobis distance between $\boldsymbol{x_{0}}$ and $\mathbf{m}_{\boldsymbol{x}_{g}}$.

Under the functional framework, the idea is to consider a similar rule but replacing the multivariate Mahalanobis distance with the functional Mahalanobis semi-distance. Consequently, assuming different means and a common
covariance operator, the new observation $\chi_{0}$ is assigned to the class $g$ if:
\begin{equation}
d_{FM}^{K}\left( \chi _{0},\widehat{\mu }_{\chi _{g}}\right) ^{2}-2\log \pi
_{g}  \label{FLBRC}
\end{equation}
is minimum. Note that the values of $\pi _{g}$ are usually fixed as the proportion of observations in the sample in the
classes. In particular, if $\pi _{1}=\cdots =\pi _{G}$, the linear Bayes classification rule reduces to the centroid classifier with the functional Mahalanobis semi-distance assuming a common covariance operator.

On the other hand, if in the multivariate case the $f_{g}$ densities are assumed to be Gaussian with different means $\mathbf{m}_{\boldsymbol{x}_{g}}$ and different covariance matrices $\mathbf{C}_{\boldsymbol{x}_{g}}$, the
Bayes rule classifies $\boldsymbol{x_{0}}$ in class $g$ if:
\begin{equation*}
d_{M}\left( \boldsymbol{x_{0}},\mathbf{m}_{\boldsymbol{x}_{g}}\right)^{2}+\log \left\vert \mathbf{C}_{\boldsymbol{x}_{g}}\right\vert -2\log \pi _{g}
\end{equation*}%
is minimum, where in this case, $d_{M}\left( \boldsymbol{x_{0}},\mathbf{m}_{\boldsymbol{x}_{g}}\right)^{2} =\left( \boldsymbol{x_{0}}-\mathbf{m}_{\boldsymbol{x}_{g}}\right) ^{\prime }\mathbf{C}_{\boldsymbol{x}_{g}}^{-1}\left( \boldsymbol{x_{0}}-\mathbf{m}_{\boldsymbol{x}_{g}}\right) $, is the squared Mahalanobis distance between $\boldsymbol{x_{0}}$ and $\mathbf{m}_{\boldsymbol{x}_{g}}$. Under the functional framework, the new observation $\chi _{0}$ is assigned to the class $G_{0}$ if:
\begin{equation}
d_{FM}^{K}(\chi ,\widehat{\mu }_{\chi _{g}})^{2}+\sum_{k=1}^{K}\log (\widehat{\lambda }_{gk})-2\log \pi _{g},  \label{FQBRC}
\end{equation}
is minimum, where $\widehat{\lambda }_{gk}$, for $k=1,\ldots ,K$ are the eigenvalues of the estimated covariance operators under class $g$, respectively, and $K$ is the number of eigenfunctions used to compute the
functional Mahalanobis semi-distances.

It is important to note that although the functional linear and quadratic classification Bayes rules in (\ref{FLBRC}) and (\ref{FQBRC}) have been derived using the functional Mahalanobis semi-distance, these methods essentially
consists in applying the multivariate linear and quadratic Bayes rules to the functional principal components scores, that are multivariate random variables. Hall et al. (2001) proposed to use the Bayes classification rule in (\ref{PP}) after estimating nonparametrically the density function of the functional principal components scores. However, these authors pointed out that a computationally less expensive method is to use the
multivariate quadratic Bayes classification rule which is essentially the rule given in (\ref{FQBRC}).

\section{Empirical results}

This section illustrates the performance of the functional classification procedures presented in Section 3 through several Monte Carlo simulations using four different scenarios and the analysis of two real datasets.

\subsection{Monte Carlo Study}

The Monte Carlo study considers four different scenarios. The first scenario consists in two Gaussian processes defined in the closed interval $I=\left[ 0,1\right] $, with different means, $\mu_{1}(t)=20t^{1.1}(1-t)$ and $\mu _{2}(t)=20t(1-t)^{1.1}$, respectively, and a common covariance operator with eigenfunctions $\psi _{k}\left( t\right) =\sqrt{2}\sin \left( \left( k-0.5\right) \pi t\right) $ and associated eigenvalues $\lambda _{k}=1/\left( \left(k-0.5\right) \pi \right) ^{2}$, for $k=1,2,\ldots $  Then, $1000$ datasets are generated composed of $n_{1}$ functions from the first process and $n_{2}$ functions from the second process such that
$n=n_{1}+n_{2}$ is the whole sample size. The generated functions are observed at $J$ equidistant points of the closed interval $I=\left[ 0,1\right] $, where $J$ is either $50$ or $100$. A Gaussian noise of variance $0.01$ is
added to each generated point. Then, once a dataset is generated in this way, the sample is split in a training sample and a test sample. The training sample is composed of $n_{10}$ functions of the first process and $n_{20}$
functions of the second process, while the test sample is composed of $n_{11}$ functions of the first process and $n_{21}$ functions of the second process such that $n_{10}+n_{11}=n_{1}$ and $n_{20}+n_{21}=n_{2}$, respectively.
In particular, two different configurations are considered. In the first one, $n=200$ with $n_{1}=n_{2}=100$ and $n_{10}=n_{20}=75$, respectively. In the second one, $n=300$ with $n_{1}=n_{2}=150$, and $n_{10}=n_{20}=120$, respectively.

The second scenario is similar to the first one but the eigenvalues of the covariance operator are given by $\lambda _{1k}=1/((k-0.5)\pi )^{2}$ and $\lambda _{2k}=2/((k-0.5)\pi )^{2}$, for $k=1,2,\ldots $, for the first and
second processes, respectively. Finally, the third and fourth scenarios are similar to the first and second ones but replacing the Gaussian process with a standardized exponential process with rate 1 and with the same mean
functions and covariance operators. The discrete trajectories are converted to functional observations using a B-splines basis of order $6$ with $20$ basis functions that are enough to fit well the data. Figure 1 shows four
datasets, once smoothing has been performed, corresponding to the four situations considered. As it can be seen in the figure, the four scenarios appear to be complicated scenarios for classification purposes.

\begin{figure}
\centering
\includegraphics[width=6in,height=3.5in]{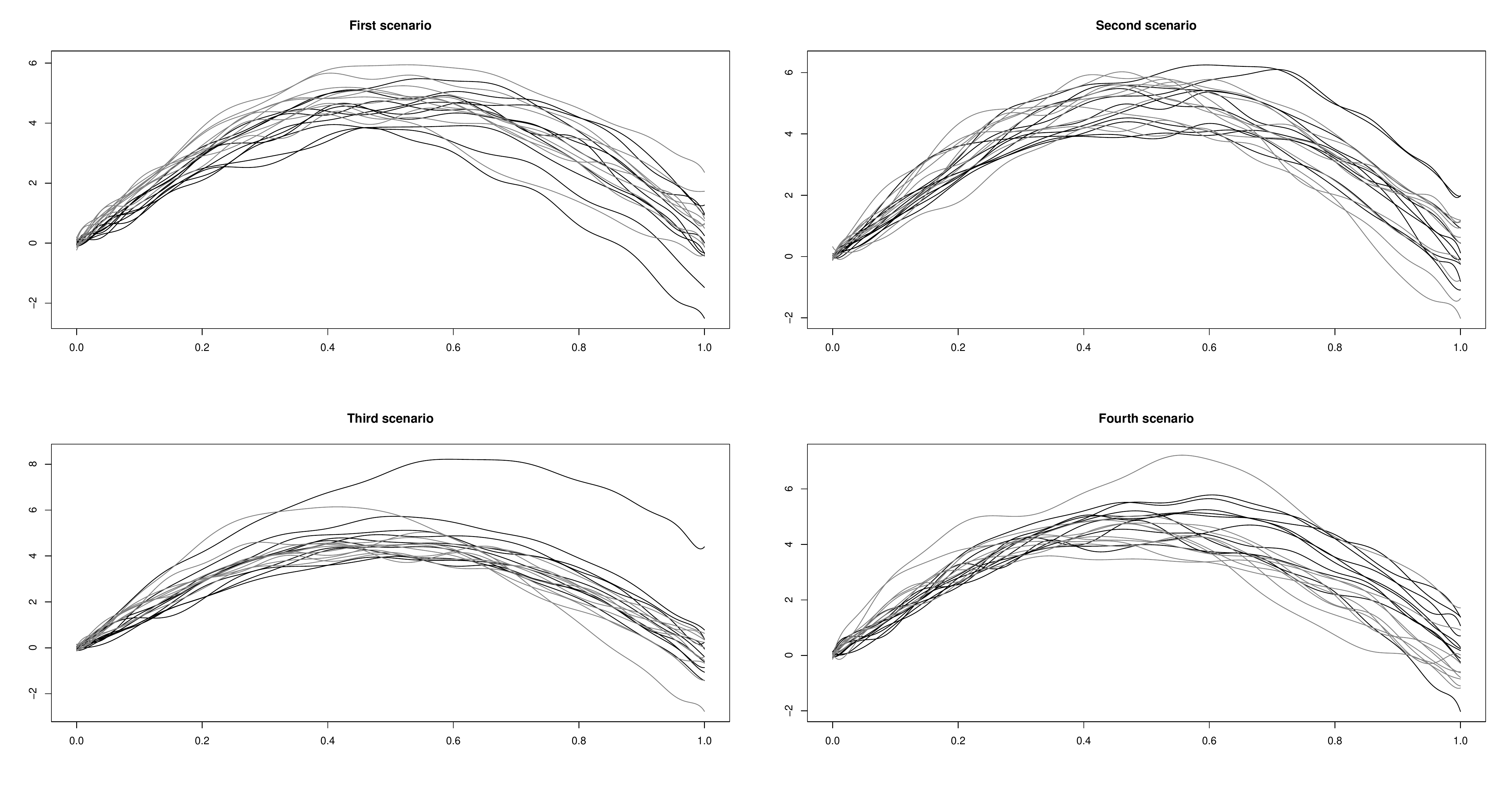}
\caption{B-splines approximations of datasets corresponding to the four experiments considered. There are $10$ functions per generated process.}
\end{figure}

For each generated dataset, the functional observations in the test sample are classified using the following procedures: (1) the kNN procedure with seven different functional distances, the $L_1$, $L_2$ and $L_{\infty}$
distances as proposed by Ba\'illo et al. (2011), the functional principal components (FPC) semi-distance assuming either a common or a different covariance operator, denoted by $FPC_C$ and $FPC_D$, respectively, and the
functional Mahalanobis (FM) semi-distance assuming either a common or a different covariance operator, denoted by $FM_C$ and $FM_D$, respectively, as proposed in Section 3; (2) the centroid procedure with eight different functional distances, the first seven as in the kNN procedure and the distance proposed by Delaigle and Hall (2012) given in (\ref{DH}) and denoted by $DH$; (3) the linear and quadratic Bayes classification rules as proposed in Section 3, denoted by $FLBCR$ and $FQBCR$, respectively; and (4) the multivariate linear and quadratic Bayes classification rules applied on the coefficients of the B-splines basis representation, denoted by $LBCR$ $Coef.$ and $QBCR$ $Coef.$, respectively. This method can be seen as a simplification of the method proposed by James and Hastie (2001) much easier to implement than the original method. The threshold values needed to compute the $FPC_C$, $FPC_D$, $FM_C$, $FM_D$ and $DH$ semi-distances and the $FLBCR$ and $FQBCR$ methods, and the maximum number of neighbors in the kNN procedures are determined using cross-validation with a maximum number of $15$ eigenfunctions and $9$ neighbors, respectively. Tables $1$, $3$, $5$ and $7$ show the proportion of correct classification of the test samples for the four scenarios. More precisely, each cell in the table displays the mean and the standard deviation (between parentheses) of the proportion of correct classifications over the $1000$ Monte Carlo samples. On the other hand, Tables $2$, $4$, $6$ and $8$ show the means and standard deviations (between parentheses) of the optimal number of principal components needed to compute the $FPC_C$, $FPC_D$, $FM_C$, $FM_D$ and $DH$ semi-distances and the $FLBCR$ and $FQBCR$ methods. In view of these tables, several comments are in order. First, in most of the cases, the kNN procedure with the $FM_C$ semi-distance attains the largest proportion of correct classifications. Second, the proportions of correct classifications for the third and fourth scenarios are larger than the corresponding proportions for the first and second scenarios suggesting that Gaussianity is not necessarily an advantage for the functional Mahalanobis semi-distance. Third, in all the situations, classification methods in conjunction with the functional Mahalanobis semi-distance have a better performance than in conjunction with any other functional distance or semi-distance or any other alternative method as the one based on the basis functions coefficients. Fourth, there is not much difference in the results in terms of number of points in the grid and sample size. Fifth, at least in these scenarios, the use of the $FPC_D$ and $FM_D$ semi-distances is not of practical advantage. Indeed, even if the generated processes have different covariance operators, the methods appear to work better assuming a common covariance operator. Sixth, note that the multivariate quadratic Bayes classification rule for the coefficients of the Basis expansion has a bad performance in all the situations. This is probably due to the large amount of parameters that is necessary to estimate. Dimension reduction as done in James and Hastie (2001) may be a solution but at the cost of increasing the complexity of the procedure. In this sense note that very simple methods provides with very good performances without a high level of sophistication. Seventh, note also that, in most of the situations, standard deviations of good classification rates linked to method based on the functional Mahalanobis semi-distance are smaller than using any other alternative. Finally, note that there is no a clear pattern relative to the number of functional principal components used with the $FPC_C$, $FPC_D$, $FM_C$, $FM_D$ and $DH$ semi-distances nor with the $FLBCR$ and $FQBCR$ methods. In summary, this limited simulation analysis appears to confirm that the functional Mahalanobis semi-distance may be a useful tool for classifying functional observations.


\begin{table}
\begin{center}
\caption{Proportion of correct classification for the first scenario}
\begin{tabular}{cccccccccccc}
\hline
$n$ & $J$ & \textbf{Method} & $L^{1}$ & $L^{2}$ & $L^{\infty }$ & $FPC_C$
& $FPC_{D}$ & $FM_C$ & $FM_{D}$ & $DH$ & $-$ \\ \hline
&  & \textbf{kNN} & $\underset{\left( .0550\right) }{.7657}$ & $\underset{%
\left( .0547\right) }{.7655}$ & $\underset{\left( .0574\right) }{.7682}$ & $%
\underset{\left( .0525\right) }{.7866}$ & $\underset{\left( .0513\right) }{%
.7871}$ & $\underset{\left( .0444\right) }{.8314}$ & $\underset{\left(
.0513\right) }{.8209}$ & $-$ & $-$ \\
&  & \textbf{Centroid} & $\underset{\left( .0870\right) }{.6710}$ & $%
\underset{\left( .0863\right) }{.6823}$ & $\underset{\left( .0854\right) }{%
.6764}$ & $\underset{\left( .0860\right) }{.6868}$ & $\underset{\left(
.0861\right) }{.6907}$ & $\underset{\left( .0490\right) }{\mathbf{.8326}}$ &
$\underset{\left( .0480\right) }{.8145}$ & $\underset{\left( .0570\right) }{%
.8017}$ & $-$ \\
200 & 50 & \textbf{FLBCR} & $-$ & $-$ & $-$ & $-$ & $-$ & $-$ & $-$ & $-$ & $%
\underset{\left( .0490\right) }{\mathbf{.8326}}$ \\
&  & \textbf{FQBCR} & $-$ & $-$ & $-$ & $-$ & $-$ & $-$ & $-$ & $-$ & $%
\underset{\left( .0480\right) }{.8145}$ \\
&  & \textbf{LBCR Coef.} & $-$ & $-$ & $-$ & $-$ & $-$ & $-$ & $-$ & $-$ & $%
\underset{\left( .0564\right) }{.8201}$ \\
&  & \textbf{QBCR Coef.} & $-$ & $-$ & $-$ & $-$ & $-$ & $-$ & $-$ & $-$ & $%
\underset{\left( .0708\right) }{.7135}$ \\ \hline
&  & \textbf{kNN} & $\underset{\left( .0584\right) }{.7700}$ & $\underset{%
\left( .0588\right) }{.7721}$ & $\underset{\left( .0553\right) }{.7744}$ & $%
\underset{\left( .0584\right) }{.7924}$ & $\underset{\left( .0570\right) }{%
.7918}$ & $\underset{\left( .0463\right) }{\mathbf{.8359}}$ & $\underset{%
\left( .0570\right) }{.8220}$ & $-$ & $-$ \\
&  & \textbf{Centroid} & $\underset{\left( .0920\right) }{.6806}$ & $%
\underset{\left( .0832\right) }{.6869}$ & $\underset{\left( .0806\right) }{%
.6837}$ & $\underset{\left( .0817\right) }{.6916}$ & $\underset{\left(
.0824\right) }{.6947}$ & $\underset{\left( .0542\right) }{.8339}$ & $%
\underset{\left( .0539\right) }{.8174}$ & $\underset{\left( .0625\right) }{%
.8061}$ & $-$ \\
200 & 100 & \textbf{FLBCR} & $-$ & $-$ & $-$ & $-$ & $-$ & $-$ & $-$ & $-$ & $%
\underset{\left( .0542\right) }{.8339}$ \\
&  & \textbf{FQBCR} & $-$ & $-$ & $-$ & $-$ & $-$ & $-$ & $-$ & $-$ & $%
\underset{\left( .0539\right) }{.8174}$ \\
&  & \textbf{LBCR Coef.} & $-$ & $-$ & $-$ & $-$ & $-$ & $-$ & $-$ & $-$ & $%
\underset{\left( .0552\right) }{.8254}$ \\
&  & \textbf{QBCR Coef.} & $-$ & $-$ & $-$ & $-$ & $-$ & $-$ & $-$ & $-$ & $%
\underset{\left( .0646\right) }{.7255}$ \\ \hline
&  & \textbf{kNN} & $\underset{\left( .0523\right) }{.7710}$ & $\underset{%
\left( .0524\right) }{.7745}$ & $\underset{\left( .0490\right) }{.7834}$ & $%
\underset{\left( .0496\right) }{.7985}$ & $\underset{\left( .0510\right) }{%
.7975}$ & $\underset{\left( .0452\right) }{.8335}$ & $\underset{\left(
.0510\right) }{.8233}$ & $-$ & $-$ \\
&  & \textbf{Centroid} & $\underset{\left( .0794\right) }{.6771}$ & $%
\underset{\left( .0752\right) }{.6853}$ & $\underset{\left( .0729\right) }{%
.6835}$ & $\underset{\left( .0761\right) }{.6897}$ & $\underset{\left(
.0762\right) }{.6915}$ & $\underset{\left( .0468\right) }{\mathbf{.8350}}$ &
$\underset{\left( .0457\right) }{.8239}$ & $\underset{\left( .0536\right) }{%
.8049}$ & $-$ \\
300 & 50 & \textbf{FLBCR} & $-$ & $-$ & $-$ & $-$ & $-$ & $-$ & $-$ & $-$ & $%
\underset{\left( .0468\right) }{\mathbf{.8350}}$ \\
&  & \textbf{FQBCR} & $-$ & $-$ & $-$ & $-$ & $-$ & $-$ & $-$ & $-$ & $%
\underset{\left( .0457\right) }{.8239}$ \\
&  & \textbf{LBCR Coef.} & $-$ & $-$ & $-$ & $-$ & $-$ & $-$ & $-$ & $-$ & $%
\underset{\left( .0488\right) }{.8325}$ \\
&  & \textbf{QBCR Coef.} & $-$ & $-$ & $-$ & $-$ & $-$ & $-$ & $-$ & $-$ & $%
\underset{\left( .0545\right) }{.7660}$ \\ \hline
&  & \textbf{kNN} & $\underset{\left( .0529\right) }{.7751}$ & $\underset{%
\left( .0521\right) }{.7766}$ & $\underset{\left( .0520\right) }{.7826}$ & $%
\underset{\left( .0523\right) }{.7948}$ & $\underset{\left( .0492\right) }{%
.7943}$ & $\underset{\left( .0450\right) }{\mathbf{.8378}}$ & $\underset{%
\left( .0492\right) }{.8225}$ & $-$ & $-$ \\
&  & \textbf{Centroid} & $\underset{\left( .0838\right) }{.6935}$ & $%
\underset{\left( .0702\right) }{.6906}$ & $\underset{\left( .0713\right) }{%
.6925}$ & $\underset{\left( .0695\right) }{.6936}$ & $\underset{\left(
.0693\right) }{.6966}$ & $\underset{\left( .0448\right) }{.8348}$ & $%
\underset{\left( .0445\right) }{.8231}$ & $\underset{\left( .0529\right) }{%
.8063}$ & $-$ \\
300 & 100 & \textbf{FLBCR} & $-$ & $-$ & $-$ & $-$ & $-$ & $-$ & $-$ & $-$ & $%
\underset{\left( .0448\right) }{.8348}$ \\
&  & \textbf{FQBCR} & $-$ & $-$ & $-$ & $-$ & $-$ & $-$ & $-$ & $-$ & $%
\underset{\left( .0445\right) }{.8231}$ \\
&  & \textbf{LBCR Coef.} & $-$ & $-$ & $-$ & $-$ & $-$ & $-$ & $-$ & $-$ & $%
\underset{\left( .0489\right) }{.8290}$ \\
&  & \textbf{QBCR Coef.} & $-$ & $-$ & $-$ & $-$ & $-$ & $-$ & $-$ & $-$ & $%
\underset{\left( .0545\right) }{.7630}$ \\ \hline
\end{tabular}
\end{center}
\end{table}

\begin{table}
\begin{center}
\caption{Means and standard deviations (between parentheses) of the optimal number of principal components needed to compute the $FPC_C$, $FPC_D$, $FM_C$, $FM_D$ and $DH$ semi-distances and the $FLBCR$ and $FQBCR$ methods}
\begin{tabular}{ccccccc}
\hline
& $FPC_C$ & $FPC_{D}$ & $FM_C$ & $FM_{D}$ & DH & $-$ \\ \hline
\textbf{kNN} & $\underset{\left( 2.87\right) }{6.36}$ & $\underset{\left(
2.97\right) }{6.53}$ & $\underset{\left( 2.90\right) }{7.48}$ & $\underset{%
\left( 2.82\right) }{7.06}$ & $-$ & $-$ \\
\textbf{Centroid} & $\underset{\left( 2.09\right) }{4.16}$ & $\underset{%
\left( 2.52\right) }{4.99}$ & $\underset{\left( 2.99\right) }{7.45}$ & $%
\underset{\left( 2.87\right) }{6.50}$ & $\underset{\left( 3.06\right) }{6.48}
$ & $-$ \\
\textbf{FLBCR} & $-$ & $-$ & $-$ & $-$ & $-$ & $\underset{\left( 2.99\right) }{%
7.45}$ \\
\textbf{FQBCR} & $-$ & $-$ & $-$ & $-$ & $-$ & $\underset{\left( 2.87\right) }{%
6.50}$ \\ \hline
\textbf{kNN} & $\underset{\left( 2.63\right) }{6.20}$ & $\underset{\left(
2.68\right) }{6.67}$ & $\underset{\left( 2.85\right) }{7.35}$ & $\underset{%
\left( 2.83\right) }{6.49}$ & $-$ & $-$ \\
\textbf{Centroid} & $\underset{\left( 2.05\right) }{4.05}$ & $\underset{%
\left( 2.34\right) }{4.66}$ & $\underset{\left( 2.93\right) }{7.36}$ & $%
\underset{\left( 2.73\right) }{6.32}$ & $\underset{\left( 3.08\right) }{6.86}
$ & $-$ \\
\textbf{FLBCR} & $-$ & $-$ & $-$ & $-$ & $-$ & $\underset{\left( 2.93\right) }{%
7.36}$ \\
\textbf{FLBCR} & $-$ & $-$ & $-$ & $-$ & $-$ & $\underset{\left( 2.73\right) }{%
6.32}$ \\ \hline
\textbf{kNN} & $\underset{\left( 2.87\right) }{6.57}$ & $\underset{\left(
2.90\right) }{6.69}$ & $\underset{\left( 2.95\right) }{7.40}$ & $\underset{%
\left( 2.83\right) }{7.21}$ & $-$ & $-$ \\
\textbf{Centroid} & $\underset{\left( 2.17\right) }{4.38}$ & $\underset{%
\left( 2.49\right) }{4.87}$ & $\underset{\left( 3.11\right) }{7.48}$ & $%
\underset{\left( 2.87\right) }{6.46}$ & $\underset{\left( 3.00\right) }{6.51}
$ & $-$ \\
\textbf{FLBCR} & $-$ & $-$ & $-$ & $-$ & $-$ & $\underset{\left( 3.11\right) }{%
7.48}$ \\
\textbf{FQBCR} & $-$ & $-$ & $-$ & $-$ & $-$ & $\underset{\left( 2.87\right) }{%
6.46}$ \\ \hline
\textbf{kNN} & $\underset{\left( 2.73\right) }{6.50}$ & $\underset{\left(
2.87\right) }{6.89}$ & $\underset{\left( 2.78\right) }{7.53}$ & $\underset{%
\left( 2.87\right) }{6.93}$ & $-$ & $-$ \\
\textbf{Centroid} & $\underset{\left( 2.09\right) }{4.39}$ & $\underset{%
\left( 2.18\right) }{4.61}$ & $\underset{\left( 2.73\right) }{7.83}$ & $%
\underset{\left( 2.86\right) }{6.67}$ & $\underset{\left( 2.97\right) }{6.94}
$ & $-$ \\
\textbf{FLBCR} & $-$ & $-$ & $-$ & $-$ & $-$ & $\underset{\left( 2.73\right) }{%
7.83}$ \\
\textbf{FQBCR} & $-$ & $-$ & $-$ & $-$ & $-$ & $\underset{\left( 2.86\right) }{%
6.67}$ \\ \hline
\end{tabular}
\end{center}
\end{table}


\begin{table}
\begin{center}
\caption{Proportion of correct classification for the second scenario}
\begin{tabular}{cccccccccccc}
\hline
$n$ & $J$ & \textbf{Method} & $L^{1}$ & $L^{2}$ & $L^{\infty }$ & $FPC_C$
& $FPC_{D}$ & $FM_C$ & $FM_{D}$ & $DH$ & $-$ \\ \hline
&  & \textbf{kNN} & $\underset{\left( .0556\right) }{.7452}$ & $\underset{%
\left( .0555\right) }{.7459}$ & $\underset{\left( .0543\right) }{.7353}$ & $%
\underset{\left( .540\right) }{.7718}$ & $\underset{\left( .0525\right) }{%
.7718}$ & $\underset{\left( .0474\right) }{\mathbf{.8055}}$ & $\underset{%
\left( .0525\right) }{.7430}$ & $-$ & $-$ \\
&  & \textbf{Centroid} & $\underset{\left( .0783\right) }{.6337}$ & $%
\underset{\left( .0774\right) }{.6415}$ & $\underset{\left( .0734\right) }{%
.6430}$ & $\underset{\left( .0774\right) }{.6469}$ & $\underset{\left(
.0783\right) }{.6497}$ & $\underset{\left( .0549\right) }{.7910}$ & $%
\underset{\left( .0554\right) }{.7130}$ & $\underset{\left( .0754\right) }{%
.7544}$ & $-$ \\
200 & 50 & \textbf{FLBCR} & $-$ & $-$ & $-$ & $-$ & $-$ & $-$ & $-$ & $-$ & $%
\underset{\left( .0549\right) }{.7910}$ \\
&  & \textbf{FQBCR} & $-$ & $-$ & $-$ & $-$ & $-$ & $-$ & $-$ & $-$ & $%
\underset{\left( .0554\right) }{.7130}$ \\
&  & \textbf{LBCR Coef.} & $-$ & $-$ & $-$ & $-$ & $-$ & $-$ & $-$ & $-$ & $%
\underset{\left( .0616\right) }{.7753}$ \\
&  & \textbf{QBCR Coef.} & $-$ & $-$ & $-$ & $-$ & $-$ & $-$ & $-$ & $-$ & $%
\underset{\left( .0502\right) }{.5817}$ \\ \hline
&  & \textbf{kNN} & $\underset{\left( .0550\right) }{.7433}$ & $\underset{%
\left( .0516\right) }{.7433}$ & $\underset{\left( .0481\right) }{.7386}$ & $%
\underset{\left( .0494\right) }{.7738}$ & $\underset{\left( .0469\right) }{%
.7747}$ & $\underset{\left( .0445\right) }{\mathbf{.8058}}$ & $\underset{%
\left( .0469\right) }{.7407}$ & $-$ & $-$ \\
&  & \textbf{Centroid} & $\underset{\left( .0869\right) }{.6350}$ & $%
\underset{\left( .0885\right) }{.6439}$ & $\underset{\left( .0813\right) }{%
.6513}$ & $\underset{\left( .0868\right) }{.6531}$ & $\underset{\left(
.0865\right) }{.6566}$ & $\underset{\left( .0528\right) }{.7928}$ & $%
\underset{\left( .0582\right) }{.7207}$ & $\underset{\left( .0637\right) }{%
.7578}$ & $-$ \\
200 & 100 & \textbf{FLBCR} & $-$ & $-$ & $-$ & $-$ & $-$ & $-$ & $-$ & $-$ & $%
\underset{\left( .0528\right) }{.7928}$ \\
&  & \textbf{FQBCR} & $-$ & $-$ & $-$ & $-$ & $-$ & $-$ & $-$ & $-$ & $%
\underset{\left( .0582\right) }{.7207}$ \\
&  & \textbf{LBCR Coef.} & $-$ & $-$ & $-$ & $-$ & $-$ & $-$ & $-$ & $-$ & $%
\underset{\left( .0581\right) }{.7805}$ \\
&  & \textbf{QBCR Coef.} & $-$ & $-$ & $-$ & $-$ & $-$ & $-$ & $-$ & $-$ & $%
\underset{\left( .0514\right) }{.5748}$ \\ \hline
&  & \textbf{kNN} & $\underset{\left( .0475\right) }{.7543}$ & $\underset{%
\left( .0455\right) }{.7550}$ & $\underset{\left( .0460\right) }{.7500}$ & $%
\underset{\left( .0438\right) }{.7833}$ & $\underset{\left( .0433\right) }{%
.7825}$ & $\underset{\left( .0397\right) }{\mathbf{.8064}}$ & $\underset{%
\left( .0433\right) }{.7359}$ & $-$ & $-$ \\
&  & \textbf{Centroid} & $\underset{\left( .0714\right) }{.6552}$ & $%
\underset{\left( .0705\right) }{.6615}$ & $\underset{\left( .0752\right) }{%
.6620}$ & $\underset{\left( .0695\right) }{.6679}$ & $\underset{\left(
.0708\right) }{.6704}$ & $\underset{\left( .0494\right) }{.7940}$ & $%
\underset{\left( .0548\right) }{.7125}$ & $\underset{\left( .0550\right) }{%
.7618}$ & $-$ \\
300 & 50 & \textbf{FLBCR} & $-$ & $-$ & $-$ & $-$ & $-$ & $-$ & $-$ & $-$ & $%
\underset{\left( .0494\right) }{.7940}$ \\
&  & \textbf{FQBCR} & $-$ & $-$ & $-$ & $-$ & $-$ & $-$ & $-$ & $-$ & $%
\underset{\left( .0548\right) }{.7125}$ \\
&  & \textbf{LBCR Coef.} & $-$ & $-$ & $-$ & $-$ & $-$ & $-$ & $-$ & $-$ & $%
\underset{\left( .0533\right) }{.7897}$ \\
&  & \textbf{QBCR Coef.} & $-$ & $-$ & $-$ & $-$ & $-$ & $-$ & $-$ & $-$ & $%
\underset{\left( .0374\right) }{.5604}$ \\ \hline
&  & \textbf{kNN} & $\underset{\left( .0490\right) }{.7538}$ & $\underset{%
\left( .0510\right) }{.7563}$ & $\underset{\left( .0490\right) }{.7516}$ & $%
\underset{\left( .0488\right) }{.7826}$ & $\underset{\left( .0492\right) }{%
.7836}$ & $\underset{\left( .0408\right) }{\mathbf{.8098}}$ & $\underset{%
\left( .0492\right) }{.7385}$ & $-$ & $-$ \\
&  & \textbf{Centroid} & $\underset{\left( .0754\right) }{.6499}$ & $%
\underset{\left( .0761\right) }{.6601}$ & $\underset{\left( .0740\right) }{%
.6651}$ & $\underset{\left( .0761\right) }{.6650}$ & $\underset{\left(
.0760\right) }{.6681}$ & $\underset{\left( .0504\right) }{.7967}$ & $%
\underset{\left( .0524\right) }{.7097}$ & $\underset{\left( .0594\right) }{%
.7650}$ & $-$ \\
300 & 100 & \textbf{FLBCR} & $-$ & $-$ & $-$ & $-$ & $-$ & $-$ & $-$ & $-$ & $%
\underset{\left( .0504\right) }{.7967}$ \\
&  & \textbf{FQBCR} & $-$ & $-$ & $-$ & $-$ & $-$ & $-$ & $-$ & $-$ & $%
\underset{\left( .0524\right) }{.7097}$ \\
&  & \textbf{LBCR Coef.} & $-$ & $-$ & $-$ & $-$ & $-$ & $-$ & $-$ & $-$ & $%
\underset{\left( .0490\right) }{.7805}$ \\
&  & \textbf{QBCR Coef.} & $-$ & $-$ & $-$ & $-$ & $-$ & $-$ & $-$ & $-$ & $%
\underset{\left( .0401\right) }{.5623}$ \\ \hline
\end{tabular}
\end{center}
\end{table}

\begin{table}
\begin{center}
\caption{Means and standard deviations (between parentheses) of the optimal number of principal components needed to compute the $FPC_C$, $FPC_D$, $FM_C$, $FM_D$ and $DH$ semi-distances for the second scenario and the $FLBCR$ and $FQBCR$ methods}
\begin{tabular}{ccccccc}
\hline
& $FPC_C$ & $FPC_{D}$ & $FM_C$ & $FM_{D}$ & $DH$ & $-$ \\ \hline
\textbf{kNN} & $\underset{\left( 2.88\right) }{6.46}$ & $\underset{\left(
2.77\right) }{6.11}$ & $\underset{\left( 2.53\right) }{6.11}$ & $\underset{%
\left( 2.69\right) }{4.94}$ & $-$ & $-$ \\
\textbf{Centroid} & $\underset{\left( 2.00\right) }{3.78}$ & $\underset{%
\left( 2.49\right) }{4.33}$ & $\underset{\left( 2.98\right) }{6.89}$ & $%
\underset{\left( 1.75\right) }{3.81}$ & $\underset{\left( 3.08\right) }{6.65}
$ & $-$ \\
\textbf{FLBCR} & $-$ & $-$ & $-$ & $-$ & $-$ & $\underset{\left( 2.98\right) }{%
6.89}$ \\
\textbf{FQBCR} & $-$ & $-$ & $-$ & $-$ & $-$ & $\underset{\left( 1.75\right) }{%
3.81}$ \\ \hline
\textbf{kNN} & $\underset{\left( 2.75\right) }{6.03}$ & $\underset{\left(
2.71\right) }{5.95}$ & $\underset{\left( 2.44\right) }{5.99}$ & $\underset{%
\left( 2.53\right) }{4.98}$ & $-$ & $-$ \\
\textbf{Centroid} & $\underset{\left( 1.89\right) }{3.77}$ & $\underset{%
\left( 2.41\right) }{4.38}$ & $\underset{\left( 2.92\right) }{6.96}$ & $%
\underset{\left( 1.67\right) }{3.62}$ & $\underset{\left( 3.16\right) }{6.52}
$ & $-$ \\
\textbf{FLBCR} & $-$ & $-$ & $-$ & $-$ & $-$ & $\underset{\left( 2.92\right) }{%
6.96}$ \\
\textbf{FQBCR} & $-$ & $-$ & $-$ & $-$ & $-$ & $\underset{\left( 1.67\right) }{%
3.62}$ \\ \hline
\textbf{kNN} & $\underset{\left( 2.67\right) }{6.40}$ & $\underset{\left(
2.66\right) }{6.33}$ & $\underset{\left( 2.37\right) }{5.95}$ & $\underset{%
\left( 2.24\right) }{4.76}$ & $-$ & $-$ \\
\textbf{Centroid} & $\underset{\left( 1.97\right) }{3.99}$ & $\underset{%
\left( 2.47\right) }{4.53}$ & $\underset{\left( 2.93\right) }{7.62}$ & $%
\underset{\left( 1.20\right) }{3.44}$ & $\underset{\left( 3.11\right) }{6.81}
$ & $-$ \\
\textbf{FLBCR} & $-$ & $-$ & $-$ & $-$ & $-$ & $\underset{\left( 2.93\right) }{%
7.62}$ \\
\textbf{FQBCR} & $-$ & $-$ & $-$ & $-$ & $-$ & $\underset{\left( 1.20\right) }{%
3.44}$ \\ \hline
\textbf{kNN} & $\underset{\left( 2.91\right) }{6.06}$ & $\underset{\left(
2.74\right) }{6.16}$ & $\underset{\left( 2.51\right) }{6.40}$ & $\underset{%
\left( 2.26\right) }{4.48}$ & $-$ & $-$ \\
\textbf{Centroid} & $\underset{\left( 1.99\right) }{4.07}$ & $\underset{%
\left( 2.55\right) }{4.82}$ & $\underset{\left( 2.90\right) }{7.41}$ & $%
\underset{\left( 1.19\right) }{3.32}$ & $\underset{\left( 3.20\right) }{7.10}
$ & $-$ \\
\textbf{FLBCR} & $-$ & $-$ & $-$ & $-$ & $-$ & $\underset{\left( 2.90\right) }{%
7.41}$ \\
\textbf{FQBCR} & $-$ & $-$ & $-$ & $-$ & $-$ & $\underset{\left( 1.19\right) }{%
3.32}$ \\ \hline
\end{tabular}
\end{center}
\end{table}


\begin{table}
\begin{center}
\caption{Proportion of correct classification for the third scenario}
\begin{tabular}{cccccccccccc}
\hline
$n$ & $J$ & \textbf{Method} & $L^{1}$ & $L^{2}$ & $L^{\infty }$ & $FPC_C$
& $FPC_{D}$ & $FM_C$ & $FM_{D}$ & $DH$ & $-$ \\ \hline
&  & \textbf{kNN} & $\underset{\left( .0498\right) }{.8344}$ & $\underset{%
\left( .0504\right) }{.8397}$ & $\underset{\left( .0439\right) }{.8408}$ & $%
\underset{\left( .0445\right) }{.8654}$ & $\underset{\left( .0437\right) }{%
.8646}$ & $\underset{\left( .0379\right) }{\mathbf{.8999}}$ & $\underset{%
\left( .0437\right) }{.8842}$ & $-$ & $-$ \\
&  & \textbf{Centroid} & $\underset{\left( .0847\right) }{.6738}$ & $%
\underset{\left( .0857\right) }{.7052}$ & $\underset{\left( .0786\right) }{%
.7131}$ & $\underset{\left( .0867\right) }{.7093}$ & $\underset{\left(
.0860\right) }{.7128}$ & $\underset{\left( .0500\right) }{.8392}$ & $%
\underset{\left( .0531\right) }{.8165}$ & $\underset{\left( .0536\right) }{%
.8099}$ & $-$ \\
200 & 50 & \textbf{FLBCR} & $-$ & $-$ & $-$ & $-$ & $-$ & $-$ & $-$ & $-$ & $%
\underset{\left( .0500\right) }{.8392}$ \\
&  & \textbf{FQBCR} & $-$ & $-$ & $-$ & $-$ & $-$ & $-$ & $-$ & $-$ & $%
\underset{\left( .0531\right) }{.8165}$ \\
&  & \textbf{LBCR Coef.} & $-$ & $-$ & $-$ & $-$ & $-$ & $-$ & $-$ & $-$ & $%
\underset{\left( .0537\right) }{.8264}$ \\
&  & \textbf{QBCR Coef.} & $-$ & $-$ & $-$ & $-$ & $-$ & $-$ & $-$ & $-$ & $%
\underset{\left( .0649\right) }{.7186}$ \\ \hline
&  & \textbf{kNN} & $\underset{\left( .0477\right) }{.8409}$ & $\underset{%
\left( .0469\right) }{.8453}$ & $\underset{\left( .0488\right) }{.8444}$ & $%
\underset{\left( .0436\right) }{.8669}$ & $\underset{\left( .0434\right) }{%
.8670}$ & $\underset{\left( .0365\right) }{\mathbf{.9050}}$ & $\underset{%
\left( .0434\right) }{.8877}$ & $-$ & $-$ \\
&  & \textbf{Centroid} & $\underset{\left( .0927\right) }{.6819}$ & $%
\underset{\left( .0944\right) }{.7060}$ & $\underset{\left( .0917\right) }{%
.7183}$ & $\underset{\left( 0.954\right) }{.7100}$ & $\underset{\left(
.0949\right) }{.7147}$ & $\underset{\left( .0512\right) }{.8464}$ & $%
\underset{\left( .0513\right) }{.8240}$ & $\underset{\left( .0574\right) }{%
.8152}$ & $-$ \\
200 & 100 & \textbf{FLBCR} & $-$ & $-$ & $-$ & $-$ & $-$ & $-$ & $-$ & $-$ & $%
\underset{\left( .0512\right) }{.8464}$ \\
&  & \textbf{FQBCR} & $-$ & $-$ & $-$ & $-$ & $-$ & $-$ & $-$ & $-$ & $%
\underset{\left( .0513\right) }{.8240}$ \\
&  & \textbf{LBCR Coef.} & $-$ & $-$ & $-$ & $-$ & $-$ & $-$ & $-$ & $-$ & $%
\underset{\left( .0574\right) }{.8342}$ \\
&  & \textbf{QBCR Coef.} & $-$ & $-$ & $-$ & $-$ & $-$ & $-$ & $-$ & $-$ & $%
\underset{\left( .0648\right) }{.7259}$ \\ \hline
&  & \textbf{kNN} & $\underset{\left( .0408\right) }{.8544}$ & $\underset{%
\left( .0409\right) }{.8596}$ & $\underset{\left( .0407\right) }{.8604}$ & $%
\underset{\left( .0368\right) }{.8821}$ & $\underset{\left( .0364\right) }{%
.8794}$ & $\underset{\left( .0294\right) }{\mathbf{.9086}}$ & $\underset{%
\left( .0364\right) }{.8984}$ & $-$ & $-$ \\
&  & \textbf{Centroid} & $\underset{\left( .0806\right) }{.6957}$ & $%
\underset{\left( .0754\right) }{.7227}$ & $\underset{\left( .0704\right) }{%
.7228}$ & $\underset{\left( .0759\right) }{.7280}$ & $\underset{\left(
.0764\right) }{.7310}$ & $\underset{\left( .0456\right) }{.8484}$ & $%
\underset{\left( .0421\right) }{.8317}$ & $\underset{\left( .0499\right) }{%
.8223}$ & $-$ \\
300 & 50 & \textbf{FLBCR} & $-$ & $-$ & $-$ & $-$ & $-$ & $-$ & $-$ & $-$ & $%
\underset{\left( .0456\right) }{.8484}$ \\
&  & \textbf{FQBCR} & $-$ & $-$ & $-$ & $-$ & $-$ & $-$ & $-$ & $-$ & $%
\underset{\left( .0421\right) }{.8317}$ \\
&  & \textbf{LBCR Coef.} & $-$ & $-$ & $-$ & $-$ & $-$ & $-$ & $-$ & $-$ & $%
\underset{\left( .0523\right) }{.8415}$ \\
&  & \textbf{QBCR Coef.} & $-$ & $-$ & $-$ & $-$ & $-$ & $-$ & $-$ & $-$ & $%
\underset{\left( .0576\right) }{.7644}$ \\ \hline
&  & \textbf{kNN} & $\underset{\left( .0405\right) }{.8570}$ & $\underset{%
\left( .0389\right) }{.8640}$ & $\underset{\left( .0433\right) }{.8621}$ & $%
\underset{\left( .0346\right) }{.8864}$ & $\underset{\left( .0354\right) }{%
.8861}$ & $\underset{\left( .0338\right) }{\mathbf{.9119}}$ & $\underset{%
\left( .0354\right) }{.9023}$ & $-$ & $-$ \\
&  & \textbf{Centroid} & $\underset{\left( .0868\right) }{.7065}$ & $%
\underset{\left( .0834\right) }{.7340}$ & $\underset{\left( .0721\right) }{%
.7363}$ & $\underset{\left( .0838\right) }{.7378}$ & $\underset{\left(
.0834\right) }{.7421}$ & $\underset{\left( .0461\right) }{.8503}$ & $%
\underset{\left( .0455\right) }{.8299}$ & $\underset{\left( .0526\right) }{%
.8245}$ & $-$ \\
300 & 100 & \textbf{FLBCR} & $-$ & $-$ & $-$ & $-$ & $-$ & $-$ & $-$ & $-$ & $%
\underset{\left( .0461\right) }{.8503}$ \\
&  & \textbf{FQBCR} & $-$ & $-$ & $-$ & $-$ & $-$ & $-$ & $-$ & $-$ & $%
\underset{\left( .0455\right) }{.8299}$ \\
&  & \textbf{LBCR Coef.} & $-$ & $-$ & $-$ & $-$ & $-$ & $-$ & $-$ & $-$ & $%
\underset{\left( .0495\right) }{.8464}$ \\
&  & \textbf{QBCR Coef.} & $-$ & $-$ & $-$ & $-$ & $-$ & $-$ & $-$ & $-$ & $%
\underset{\left( .0565\right) }{.7623}$ \\ \hline
\end{tabular}
\end{center}
\end{table}

\begin{table}
\begin{center}
\caption{Means and standard deviations (between parentheses) of the optimal number of principal components needed to compute the $FPC_C$, $FPC_D$, $FM_C$, $FM_D$ and $DH$ semi-distances for the third scenario and the $FLBCR$ and $FQBCR$ methods}
\begin{tabular}{ccccccc}
\hline
& $FPC_C$ & $FPC_{D}$ & $FM_C$ & $FM_{D}$ & $DH$ & $-$ \\ \hline
\textbf{kNN} & $\underset{\left( 2.84\right) }{5.71}$ & $\underset{\left(
2.74\right) }{6.04}$ & $\underset{\left( 2.59\right) }{5.35}$ & $\underset{%
\left( 2.45\right) }{5.00}$ & $-$ & $-$ \\
\textbf{Centroid} & $\underset{\left( 2.24\right) }{4.25}$ & $\underset{%
\left( 2.26\right) }{4.73}$ & $\underset{\left( 2.77\right) }{7.20}$ & $%
\underset{\left( 2.70\right) }{6.09}$ & $\underset{\left( 2.74\right) }{6.42}
$ & $-$ \\
\textbf{FLBCR} & $-$ & $-$ & $-$ & $-$ & $-$ & $\underset{\left( 2.77\right) }{%
7.20}$ \\
\textbf{FQBCR} & $-$ & $-$ & $-$ & $-$ & $-$ & $\underset{\left( 2.70\right) }{%
6.09}$ \\ \hline
\textbf{kNN} & $\underset{\left( 2.71\right) }{5.62}$ & $\underset{\left(
2.77\right) }{5.98}$ & $\underset{\left( 2.56\right) }{5.14}$ & $\underset{%
\left( 2.39\right) }{4.92}$ & $-$ & $-$ \\
\textbf{Centroid} & $\underset{\left( 1.98\right) }{3.91}$ & $\underset{%
\left( 2.43\right) }{4.66}$ & $\underset{\left( 2.76\right) }{6.68}$ & $%
\underset{\left( 2.60\right) }{5.67}$ & $\underset{\left( 3.00\right) }{6.46}
$ & $-$ \\
\textbf{FLBCR} & $-$ & $-$ & $-$ & $-$ & $-$ & $\underset{\left( 2.76\right) }{%
6.68}$ \\
\textbf{FQBCR} & $-$ & $-$ & $-$ & $-$ & $-$ & $\underset{\left( 2.60\right) }{%
5.67}$ \\ \hline
\textbf{kNN} & $\underset{\left( 2.84\right) }{6.10}$ & $\underset{\left(
2.89\right) }{6.29}$ & $\underset{\left( 2.37\right) }{4.93}$ & $\underset{%
\left( 2.26\right) }{4.78}$ & $-$ & $-$ \\
\textbf{Centroid} & $\underset{\left( 2.20\right) }{4.22}$ & $\underset{%
\left( 2.55\right) }{4.84}$ & $\underset{\left( 2.81\right) }{7.14}$ & $%
\underset{\left( 2.79\right) }{6.04}$ & $\underset{\left( 2.92\right) }{7.03}
$ & $-$ \\
\textbf{FLBCR} & $-$ & $-$ & $-$ & $-$ & $-$ & $\underset{\left( 2.81\right) }{%
7.14}$ \\
\textbf{FQBCR} & $-$ & $-$ & $-$ & $-$ & $-$ & $\underset{\left( 2.79\right) }{%
6.04}$ \\ \hline
\textbf{kNN} & $\underset{\left( 2.75\right) }{5.91}$ & $\underset{\left(
2.77\right) }{6.52}$ & $\underset{\left( 2.21\right) }{4.58}$ & $\underset{%
\left( 1.94\right) }{4.52}$ & $-$ & $-$ \\
\textbf{Centroid} & $\underset{\left( 1.99\right) }{4.34}$ & $\underset{%
\left( 2.40\right) }{5.08}$ & $\underset{\left( 2.97\right) }{7.23}$ & $%
\underset{\left( 2.85\right) }{6.17}$ & $\underset{\left( 2.91\right) }{6.75}
$ & $-$ \\
\textbf{FLBCR} & $-$ & $-$ & $-$ & $-$ & $-$ & $\underset{\left( 2.97\right) }{%
7.23}$ \\
\textbf{FQBCR} & $-$ & $-$ & $-$ & $-$ & $-$ & $\underset{\left( 2.85\right) }{%
6.17}$ \\ \hline
\end{tabular}
\end{center}
\end{table}


\begin{table}
\begin{center}
\caption{Proportion of correct classification for the fourth scenario}
\begin{tabular}{cccccccccccc}
\hline
$n$ & $J$ & \textbf{Method} & $L^{1}$ & $L^{2}$ & $L^{\infty }$ & $FPC_C$
& $FPC_{D}$ & $FM_C$ & $FM_{D}$ & $DH$ & $-$ \\ \hline
&  & \textbf{kNN} & $\underset{\left( .0489\right) }{.8645}$ & $\underset{%
\left( .0485\right) }{.8647}$ & $\underset{\left( .0461\right) }{.8517}$ & $%
\underset{\left( .0464\right) }{.8851}$ & $\underset{\left( .0457\right) }{%
.8830}$ & $\underset{\left( .0349\right) }{\mathbf{.9212}}$ & $\underset{%
\left( .0457\right) }{.8619}$ & $-$ & $-$ \\
&  & \textbf{Centroid} & $\underset{\left( .0849\right) }{.7328}$ & $%
\underset{\left( .0892\right) }{.7213}$ & $\underset{\left( .0903\right) }{%
.6923}$ & $\underset{\left( .0897\right) }{.7234}$ & $\underset{\left(
.0887\right) }{.7262}$ & $\underset{\left( .0424\right) }{.8939}$ & $%
\underset{\left( .0478\right) }{.8386}$ & $\underset{\left( .0447\right) }{%
.8679}$ & $-$ \\
200 & 50 & \textbf{FLBCR} & $-$ & $-$ & $-$ & $-$ & $-$ & $-$ & $-$ & $-$ & $%
\underset{\left( .0424\right) }{.8939}$ \\
&  & \textbf{FQBCR} & $-$ & $-$ & $-$ & $-$ & $-$ & $-$ & $-$ & $-$ & $%
\underset{\left( .0478\right) }{.8386}$ \\
&  & \textbf{LBCR Coef.} & $-$ & $-$ & $-$ & $-$ & $-$ & $-$ & $-$ & $-$ & $%
\underset{\left( .0450\right) }{.8913}$ \\
&  & \textbf{QBCR Coef.} & $-$ & $-$ & $-$ & $-$ & $-$ & $-$ & $-$ & $-$ & $%
\underset{\left( .0668\right) }{.7173}$ \\ \hline
&  & \textbf{kNN} & $\underset{\left( .0465\right) }{.8669}$ & $\underset{%
\left( .0447\right) }{.8712}$ & $\underset{\left( .0490\right) }{.8576}$ & $%
\underset{\left( .0414\right) }{.8874}$ & $\underset{\left( .0425\right) }{%
.8872}$ & $\underset{\left( .0330\right) }{\mathbf{.9223}}$ & $\underset{%
\left( .0425\right) }{.8599}$ & $-$ & $-$ \\
&  & \textbf{Centroid} & $\underset{\left( .0846\right) }{.7345}$ & $%
\underset{\left( .0856\right) }{.7291}$ & $\underset{\left( .0861\right) }{%
.6969}$ & $\underset{\left( .0861\right) }{.7316}$ & $\underset{\left(
.0859\right) }{.7337}$ & $\underset{\left( .0398\right) }{.8930}$ & $%
\underset{\left( .0456\right) }{.8335}$ & $\underset{\left( .0455\right) }{%
.8681}$ & $-$ \\
200 & 100 & \textbf{FLBCR} & $-$ & $-$ & $-$ & $-$ & $-$ & $-$ & $-$ & $-$ & $%
\underset{\left( .0398\right) }{.8930}$ \\
&  & \textbf{FQBCR} & $-$ & $-$ & $-$ & $-$ & $-$ & $-$ & $-$ & $-$ & $%
\underset{\left( .0456\right) }{.8335}$ \\
&  & \textbf{LBCR Coef.} & $-$ & $-$ & $-$ & $-$ & $-$ & $-$ & $-$ & $-$ & $%
\underset{\left( .0412\right) }{.8912}$ \\
&  & \textbf{QBCR Coef.} & $-$ & $-$ & $-$ & $-$ & $-$ & $-$ & $-$ & $-$ & $%
\underset{\left( .0632\right) }{.7077}$ \\ \hline
&  & \textbf{kNN} & $\underset{\left( .0407\right) }{.8803}$ & $\underset{%
\left( .0419\right) }{.8854}$ & $\underset{\left( .0430\right) }{.8713}$ & $%
\underset{\left( .0383\right) }{.9006}$ & $\underset{\left( .0387\right) }{%
.8987}$ & $\underset{\left( .0316\right) }{\mathbf{.9265}}$ & $\underset{%
\left( .0387\right) }{.8670}$ & $-$ & $-$ \\
&  & \textbf{Centroid} & $\underset{\left( .0688\right) }{.7303}$ & $%
\underset{\left( .0747\right) }{.7250}$ & $\underset{\left( .0800\right) }{%
.7040}$ & $\underset{\left( .0749\right) }{.7268}$ & $\underset{\left(
.0741\right) }{.7293}$ & $\underset{\left( .0379\right) }{.8947}$ & $%
\underset{\left( .0445\right) }{.8340}$ & $\underset{\left( .0479\right) }{%
.8647}$ & $-$ \\
300 & 50 & \textbf{FLBCR} & $-$ & $-$ & $-$ & $-$ & $-$ & $-$ & $-$ & $-$ & $%
\underset{\left( .0379\right) }{.8947}$ \\
&  & \textbf{FQBCR} & $-$ & $-$ & $-$ & $-$ & $-$ & $-$ & $-$ & $-$ & $%
\underset{\left( .0445\right) }{.8340}$ \\
&  & \textbf{LBCR Coef.} & $-$ & $-$ & $-$ & $-$ & $-$ & $-$ & $-$ & $-$ & $%
\underset{\left( .0389\right) }{.9044}$ \\
&  & \textbf{QBCR Coef.} & $-$ & $-$ & $-$ & $-$ & $-$ & $-$ & $-$ & $-$ & $%
\underset{\left( .0589\right) }{.7270}$ \\ \hline
&  & \textbf{kNN} & $\underset{\left( .0420\right) }{.8795}$ & $\underset{%
\left( .0442\right) }{.8797}$ & $\underset{\left( .0427\right) }{.8632}$ & $%
\underset{\left( .0386\right) }{.8974}$ & $\underset{\left( .0392\right) }{%
.8960}$ & $\underset{\left( .0318\right) }{\mathbf{.9279}}$ & $\underset{%
\left( 0.392\right) }{.8710}$ & $-$ & $-$ \\
&  & \textbf{Centroid} & $\underset{\left( .0717\right) }{.7346}$ & $%
\underset{\left( .0797\right) }{.7289}$ & $\underset{\left( .0818\right) }{%
.7015}$ & $\underset{\left( .0795\right) }{.7314}$ & $\underset{\left(
.0796\right) }{.7330}$ & $\underset{\left( .0384\right) }{.8936}$ & $%
\underset{\left( .0441\right) }{.8324}$ & $\underset{\left( .0455\right) }{%
.8649}$ & $-$ \\
300 & 100 & \textbf{FLBCR} & $-$ & $-$ & $-$ & $-$ & $-$ & $-$ & $-$ & $-$ & $%
\underset{\left( .0384\right) }{.8936}$ \\
&  & \textbf{FQBCR} & $-$ & $-$ & $-$ & $-$ & $-$ & $-$ & $-$ & $-$ & $%
\underset{\left( .0441\right) }{.8324}$ \\
&  & \textbf{LBCR Coef.} & $-$ & $-$ & $-$ & $-$ & $-$ & $-$ & $-$ & $-$ & $%
\underset{\left( .0402\right) }{.9010}$ \\
&  & \textbf{QBCR Coef.} & $-$ & $-$ & $-$ & $-$ & $-$ & $-$ & $-$ & $-$ & $%
\underset{\left( .0613\right) }{.7215}$ \\ \hline
\end{tabular}
\end{center}
\end{table}

\begin{table}
\begin{center}
\caption{Means and standard deviations (between parentheses) of the optimal number of principal components needed to compute the $FPC_C$, $FPC_D$, $FM_C$, $FM_D$ and $DH$ semi-distances for the fourth scenario and the $FLBCR$ and $FQBCR$ methods}
\begin{tabular}{ccccccc}
\hline
& $FPC_C$ & $FPC_{D}$ & $FM_C$ & $FM_{D}$ & $DH$ & $-$ \\ \hline
\textbf{kNN} & $\underset{\left( 2.99\right) }{6.41}$ & $\underset{\left(
2.83\right) }{6.51}$ & $\underset{\left( 3.00\right) }{6.49}$ & $\underset{%
\left( 1.56\right) }{3.72}$ & $-$ & $-$ \\
\textbf{Centroid} & $\underset{\left( 2.51\right) }{4.87}$ & $\underset{%
\left( 2.72\right) }{5.23}$ & $\underset{\left( 2.65\right) }{8.26}$ & $%
\underset{\left( 2.49\right) }{5.06}$ & $\underset{\left( 2.77\right) }{7.54}
$ & $-$ \\
\textbf{FLBCR} & $-$ & $-$ & $-$ & $-$ & $-$ & $\underset{\left( 2.65\right) }{%
8.26}$ \\
\textbf{FQBCR} & $-$ & $-$ & $-$ & $-$ & $-$ & $\underset{\left( 2.49\right) }{%
5.06}$ \\ \hline
\textbf{kNN} & $\underset{\left( 2.83\right) }{6.52}$ & $\underset{\left(
2.79\right) }{7.02}$ & $\underset{\left( 2.94\right) }{6.47}$ & $\underset{%
\left( 1.73\right) }{3.71}$ & $-$ & $-$ \\
\textbf{Centroid} & $\underset{\left( 2.53\right) }{5.07}$ & $\underset{%
\left( 2.65\right) }{5.51}$ & $\underset{\left( 2.67\right) }{8.32}$ & $%
\underset{\left( 2.38\right) }{4.62}$ & $\underset{\left( 2.89\right) }{7.77}
$ & $-$ \\
\textbf{FLBCR} & $-$ & $-$ & $-$ & $-$ & $-$ & $\underset{\left( 2.67\right) }{%
8.32}$ \\
\textbf{FQBCR} & $-$ & $-$ & $-$ & $-$ & $-$ & $\underset{\left( 2.38\right) }{%
4.62}$ \\ \hline
\textbf{kNN} & $\underset{\left( 2.89\right) }{6.96}$ & $\underset{\left(
2.84\right) }{7.11}$ & $\underset{\left( 2.90\right) }{6.16}$ & $\underset{%
\left( 1.27\right) }{3.50}$ & $-$ & $-$ \\
\textbf{Centroid} & $\underset{\left( 2.39\right) }{5.00}$ & $\underset{%
\left( 2.31\right) }{5.28}$ & $\underset{\left( 2.64\right) }{8.51}$ & $%
\underset{\left( 2.14\right) }{4.73}$ & $\underset{\left( 2.90\right) }{7.51}
$ & $-$ \\
\textbf{FLBCR} & $-$ & $-$ & $-$ & $-$ & $-$ & $\underset{\left( 2.64\right) }{%
8.51}$ \\
\textbf{FQBCR} & $-$ & $-$ & $-$ & $-$ & $-$ & $\underset{\left( 2.14\right) }{%
4.73}$ \\ \hline
\textbf{kNN} & $\underset{\left( 2.72\right) }{6.64}$ & $\underset{\left(
2.85\right) }{7.04}$ & $\underset{\left( 2.94\right) }{6.42}$ & $\underset{%
\left( 1.09\right) }{3.37}$ & $-$ & $-$ \\
\textbf{Centroid} & $\underset{\left( 2.20\right) }{4.90}$ & $\underset{%
\left( 2.67\right) }{5.55}$ & $\underset{\left( 2.55\right) }{8.42}$ & $%
\underset{\left( 2.25\right) }{4.96}$ & $\underset{\left( 2.89\right) }{7.30}
$ & $-$ \\
\textbf{FLBCR} & $-$ & $-$ & $-$ & $-$ & $-$ & $\underset{\left( 2.55\right) }{%
8.42}$ \\
\textbf{FQBCR} & $-$ & $-$ & $-$ & $-$ & $-$ & $\underset{\left( 2.25\right) }{%
4.96}$ \\ \hline
\end{tabular}
\end{center}
\end{table}

\subsection{Real data study: Tecator dataset}

Next, the classification procedures are applied to the Tecator dataset previously considered by Ferraty and Vieu (2003), Rossi and Villa (2006), Li and Yu (2008), Alonso et al. (2012) and Martin-Barragan et al. (2013), among others. The dataset that consists of $215$ near-infrared absorbance spectra of meat samples, recorded on a Tecator Infracted Food Analyzer is available at http://lib.stat.cmu.edu/datasets/tecator. The absorbance of a meat sample is a function given by $log_{10}\left(I_0/I\right)$ where $I_0$ and $I$ are, respectively, the intensity of the light before and after passing through of the meat sample. Each observation consist of a $100$-channel absorbance spectrum in the wavelength range 850-1050 nm, contents of moisture (water), fat and protein. Therefore, the recorded absorbance can be seen as a discretized version of the continuous process. The classification problem here is to separate meat samples with a high fat content (more than $20\%$) from samples with low fat content (less than $20\%$) based on the absorbance. Among the $215$ samples, $77$ have high fat content and $138$ have low fat content. Previous analysis of this dataset have suggested that classification of the second order derivatives of the observed functions produces lower misclassification rates. Therefore, the analysis of the original data and their second order derivatives are carried out. In both cases, the discrete observations are converted to functional observations using a B-splines basis of order $6$ with $20$ and $40$ basis functions, respectively, that are enough to fit well the data. Figure $2$ shows the sample of these $100$-channel absorbance spectrum and their second derivatives after smoothing.

\begin{figure}
\centering
\includegraphics[width=6in,height=3.5in]{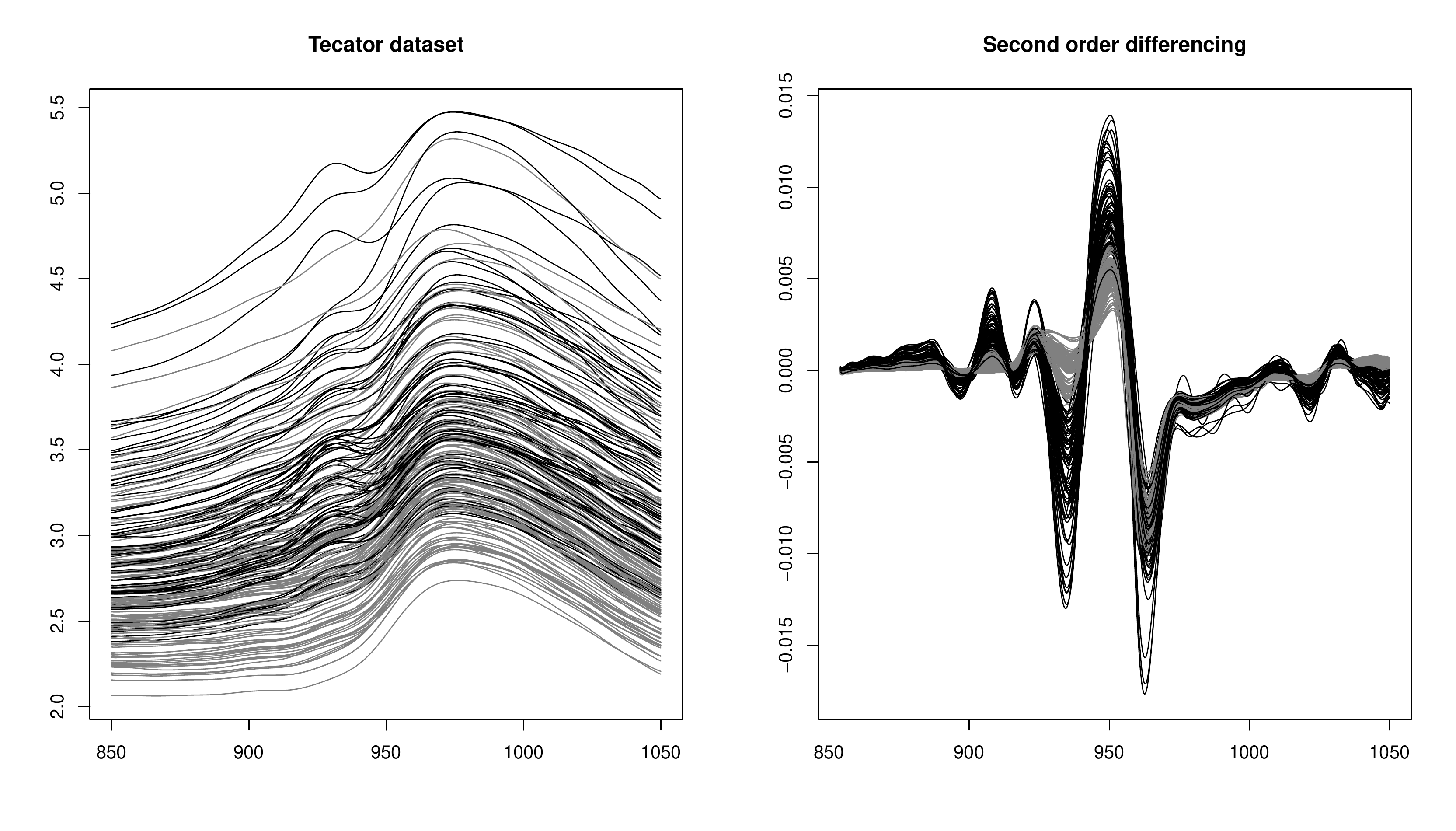}
\caption{Right: Original observations of the Tecator dataset. Left: Second order derivatives of the Tecator dataset. High fat content in black and low fat content in gray}
\end{figure}

In order to evaluate the performance of the functional classification methods given before, $1000$ training samples are considered composed by $58$ and $104$ randomly chosen functions of meat with high fat content and low fat content, respectively. For each training sample, it is associated a test sample composed by the remaining $19$ and $34$ functions of meat with high fat content and low fat content, respectively. The classification results are shown in Tables $9$ and $11$ that show the mean and the standard deviation (between parentheses) of the proportion of correct classifications obtained via cross-validation for the two cases. As in the simulation study, the threshold values needed to compute the $FPC_C$, $FPC_D$, $FM_C$, $FM_D$ and $DH$ semi-distances and the $FLBCR$ and $FQBCR$ methods, and the maximum number of neighbors in the kNN procedures are determined using cross-validation with a maximum number of $15$ eigenfunctions and $9$ neighbors, respectively. In both cases, the kNN procedure with the $FM_C$ semi-distance is the winner. The highest proportions of correct classification for the Tecator dataset and the second order derivatives are $0.9835$ and $0.9918$, respectively, suggesting that it is not necessary to use the second order derivatives of the Tecator data to obtain almost perfect classification. Note that using a similar experiment, Rossi and Villa (2006) obtained good classification rates of $0.9672$ and $0.9740$ for the original and second order derivatives with SVMs, respectively, Li and Yu (2008) obtained good classification rates of $0.9602$ and $0.9891$ for the original and second order derivatives with a segmentation approach, respectively, Alonso et al. (2012) obtained good classification rates of $0.9798$ and $0.9768$, respectively, with two methods that takes into account the original, the first and the second order derivatives, and, finally, Martin-Barragan et al. (2013) obtained a good classification rate of $0.9891$ with SVMs. Note that all of the previous approaches are more sophisticated than the ones taken here.

On the other hand, Tables $10$ and $12$ show the means and the standard deviation (between parentheses) of the number of principal components needed to calculate the $FPC_C$, $FPC_D$, $FM_C$, $FM_D$ and $DH$ semi-distances and the $FLBCR$ and $FQBCR$ methods for the original dataset and their second order derivatives. The mean numbers of functional principal components used with the functional Mahalanobis semi-distance are slightly larger than the corresponding to the functional principal components and Delaigle and Hall semi-distances if the original dataset is used but are sometimes smaller for their second order derivatives. Therefore, apparently there is not a general rule regarding the number of principal components used.

\begin{table}
\begin{center}
\caption{Proportion of correct classification for the Tecator dataset}
\begin{tabular}{cccccccccc}
\hline
& $L^{1}$ & $L^{2}$ & $L^{\infty }$ & $FPC_C$ & $FPC_{D}$ & $FM_C$ & $FM_{D}$ & $DH$ & $-$ \\ \hline
\textbf{kNN} & $\underset{\left( .0368\right) }{.7904}$ & $\underset{\left(
.0371\right) }{.8108}$ & $\underset{\left( .0342\right) }{.8602}$ & $%
\underset{\left( .0364\right) }{.8144}$ & $\underset{\left( .0363\right) }{%
.8135}$ & $\underset{\left( .0114\right) }{\mathbf{.9835}}$ & $\underset{%
\left( .0363\right) }{.9714}$ & $-$ & $-$ \\
\textbf{Centroid} & $\underset{\left( .0343\right) }{.6784}$ & $\underset{%
\left( .0347\right) }{.6812}$ & $\underset{\left( .0346\right) }{.6957}$ & $%
\underset{\left( .0347\right) }{.6813}$ & $\underset{\left( .0348\right) }{%
.6813}$ & $\underset{\left( .0173\right) }{.9630}$ & $\underset{\left(
.0218\right) }{.9521}$ & $\underset{\left( .0322\right) }{.9479}$ & $-$ \\
\textbf{FLBCR} & $-$ & $-$ & $-$ & $-$ & $-$ & $-$ & $-$ & $-$ & $\underset{\left( .0196\right)
}{.9517}$ \\
\textbf{FQBCR} & $-$ & $-$ & $-$ & $-$ & $-$ & $-$ & $-$ & $-$ & $\underset{\left( .0172\right) }{.9671%
}$ \\
\textbf{LBCR Coef.} & $-$ & $-$ & $-$ & $-$ & $-$ & $-$ & $-$ & $-$ & $\underset{\left( .0244\right) }%
{.9244}$ \\
\textbf{QBCR Coef.} & $-$ & $-$ & $-$ & $-$ & $-$ & $-$ & $-$ & $-$ & $\underset{\left(
.0325\right) }{.8958}$ \\ \hline
\end{tabular}
\end{center}
\end{table}

\begin{table}
\begin{center}
\caption{Means and standard deviations of the number of principal components used by the $FPC_C$, $FPC_D$, $FM_C$, $FM_D$ and $DH$ semi-distances and the $FLBCR$ and $FQBCR$ methods for the Tecator dataset}
\begin{tabular}{ccccccc}
\hline
& $FPC_C$ & $FPC_{D}$ & $FM_C$ & $FM_{D}$ & $DH$ & $-$ \\ \hline
\textbf{kNN} & $\underset{\left( .073\right) }{4.19}$ & $\underset{\left(
1.54\right) }{4.89}$ & $\underset{\left( 1.01\right) }{4.86}$ & $\underset{%
\left( 1.18\right) }{5.12}$ & $-$ & $-$ \\
\textbf{Centroid} & $\underset{\left( 0.98\right) }{1.48}$ & $\underset{%
\left( 1.09\right) }{1.52}$ & $\underset{\left( 0.94\right) }{4.82}$ & $%
\underset{\left( 1.24\right) }{5.20}$ & $\underset{\left( 1.47\right) }{5.05}
$ & $-$ \\
\textbf{FLBCR} & $-$ & $-$ & $-$ & $-$ & $-$ & $\underset{\left( 1.10\right) }{%
5.01}$ \\
\textbf{FQBCR} & $-$ & $-$ & $-$ & $-$ & $-$ & $\underset{\left( 1.10\right) }{%
5.15}$ \\ \hline
\end{tabular}
\end{center}
\end{table}

\begin{table}
\begin{center}
\caption{Proportion of correct classification for the second order differences of the Tecator dataset}
\begin{tabular}{cccccccccc}
\hline
& $L^{1}$ & $L^{2}$ & $L^{\infty }$ & $FPC_C$ & $FPC_{D}$ & $FM_C$ & $FM_{D}$ & $DH$ & $-$ \\ \hline
\textbf{kNN} & $\underset{\left( .0091\right) }{.9885}$ & $\underset{\left(
.0099\right) }{.9852}$ & $\underset{\left( .0109\right) }{.9814}$ & $%
\underset{\left( .0080\right) }{.9901}$ & $\underset{\left( .0094\right) }{%
.9870}$ & $\underset{\left( .0076\right) }{\mathbf{.9918}}$ & $\underset{%
\left( .0094\right) }{.9664}$ & $-$ & $-$ \\
\textbf{Centroid} & $\underset{\left( .0200\right) }{.9629}$ & $\underset{%
\left( .0210\right) }{.9608}$ & $\underset{\left( .0217\right) }{.9546}$ & $%
\underset{\left( .0190\right) }{.9651}$ & $\underset{\left( .0206\right) }{%
.9617}$ & $\underset{\left( .0180\right) }{.9678}$ & $\underset{\left(
.0253\right) }{.9372}$ & $\underset{\left( .0201\right) }{.9630}$ &  $-$ \\
\textbf{FLBCR} & $-$ & $-$ & $-$ & $-$ & $-$ & $-$ & $-$ & $-$ & $\underset{\left( .0195\right)
}{.9533}$ \\
\textbf{FQBCR} & $-$ & $-$ & $-$ & $-$ & $-$ & $-$ & $-$ & $-$ & $\underset{\left( .0190\right) }{.9555%
}$ \\
\textbf{LBCR Coef.} & $-$ & $-$ & $-$ & $-$ & $-$ & $-$ & $-$ & $-$ & $\underset{\left( .0261\right) }%
{.9218}$ \\
\textbf{QBCR Coef.} & $-$ & $-$ & $-$ & $-$ & $-$ & $-$ & $-$ & $-$ & $\underset{\left(
.0581\right) }{.7220}$ \\ \hline
\end{tabular}
\end{center}
\end{table}

\begin{table}
\begin{center}
\caption{Means and standard deviations of the number of principal components used by the $FPC_C$, $FPC_D$, $FM_C$, $FM_D$ and $DH$ semi-distances and the $FLBCR$ and $FQBCR$ methods for the second order derivatives of the Tecator dataset}
\begin{tabular}{ccccccc}
\hline
& $FPC_C$ & $FPC_{D}$ & $FM_C$ & $FM_{D}$ & $DH$ & $-$ \\ \hline
\textbf{kNN} & $\underset{\left( 0.62\right) }{2.22}$ & $\underset{\left(
1.99\right) }{3.78}$ & $\underset{\left( 1.83\right) }{2.67}$ & $\underset{%
\left( 1.00\right) }{2.05}$ & $-$ & $-$ \\
\textbf{Centroid} & $\underset{\left( 0.63\right) }{1.63}$ & $\underset{%
\left( 1.73\right) }{3.35}$ & $\underset{\left( 2.96\right) }{2.99}$ & $%
\underset{\left( 1.43\right) }{1.66}$ & $\underset{\left( 3.47\right) }{3.59}
$ & $-$ \\
\textbf{FLBCR} & $-$ & $-$ & $-$ & $-$ & $-$ & $\underset{\left( 3.85\right) }{%
4.24}$ \\
\textbf{FQBCR} & $-$ & $-$ & $-$ & $-$ & $-$ & $\underset{\left( 1.29\right) }{%
2.00}$ \\ \hline
\end{tabular}
\end{center}
\end{table}

\subsection{Real data study: Phoneme dataset}

Finally, the classification procedures are applied to the Phoneme dataset described in Ferraty and Vieu (2006) and available at http://www.math.univ-toulouse.fr/staph/npfda/npfda-datasets.html. The dataset contains log-periodograms corresponding to recordings of speakers of $32$ ms duration. Here, two populations are considered corresponding to the phonemes ``aa'' as the vowel in ``dark'' and ``ao'' as the first vowel in ``water'', such
that each speech frame is represented by $400$ samples at a $16$-kHz sampling rate where only the first $150$ frequencies from each subject are retained. Therefore, the data consists of $800$ log-periodograms of length $150$,
with known class phoneme membership. The classification problem here is to separate the two phonemes. The discrete observations are converted to functional observations using a B-splines basis of order $6$ with $40$ basis
functions, respectively, that are enough to fit well the data. Figure $3$ shows the sample of log-periodograms. The figure confirms that it is difficult to distinguish the log-periodograms from one another.

\begin{figure}
\centering
\includegraphics[width=6in,height=3.5in]{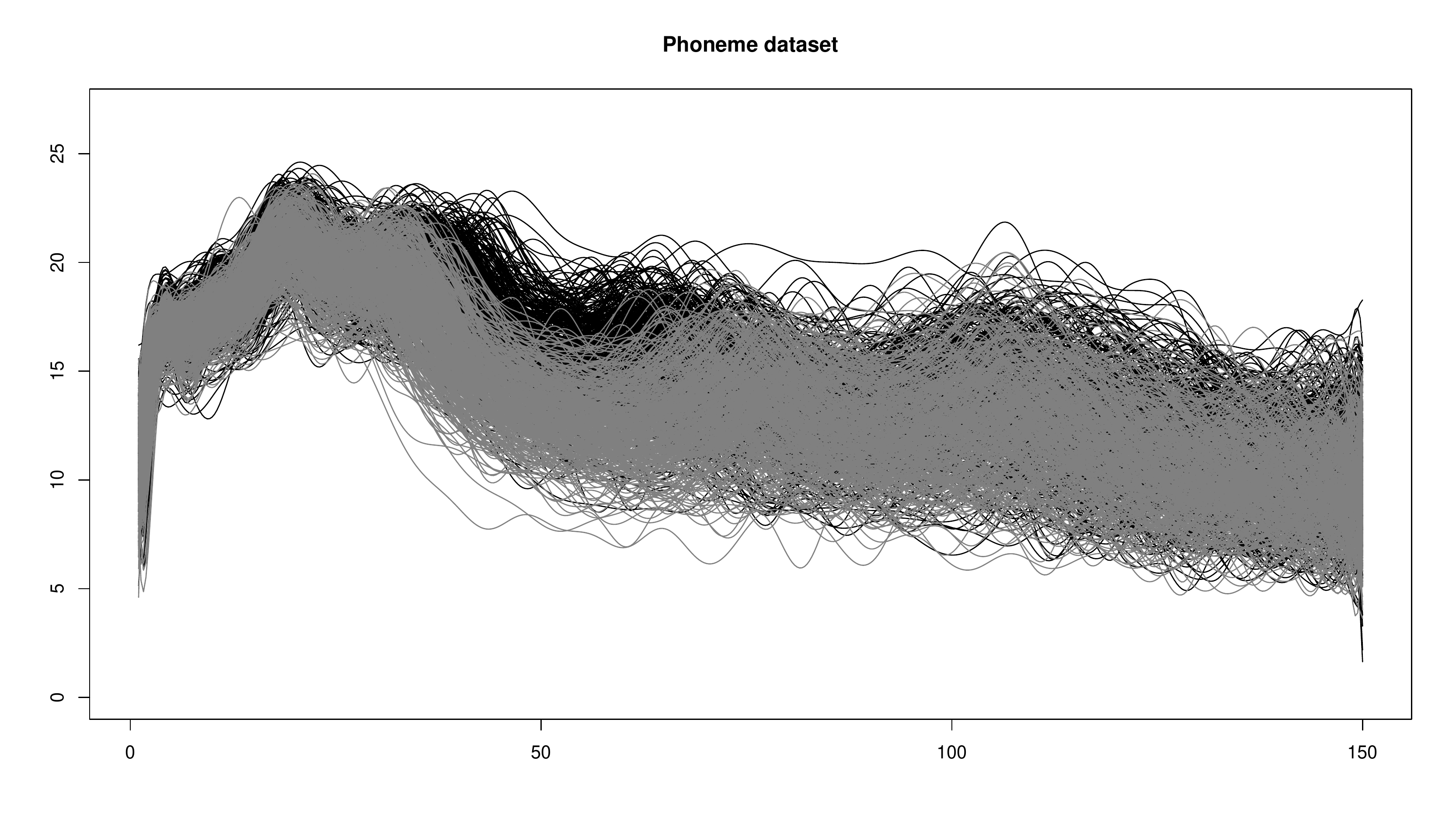}
\caption{Phoneme dataset. Log-periodograms for ``aa'' in black and log-periodograms for ``ao'' in gray. Note that the log-periodograms in gray hide most of the log-periodograms in black}
\end{figure}

As in the previous example, $1000$ training samples are considered composed by $300$ randomly chosen log-periodograms of both vowels. For each training sample, it is associated a test sample composed by the remaining $200$ log-periodograms, $100$ per vowel, respectively. The classification results are shown in Table $13$ that shows the mean and the standard deviation (between parentheses) of the proportion of correct classifications obtained via cross-validation. As in the simulation study and the previous example, the threshold values needed to compute the $FPC_C$, $FPC_D$, $FM_C$, $FM_D$ and $DH$ semi-distances and the $FLBCR$ and $FQBCR$ methods, and the maximum number of neighbors in the kNN procedures are determined using cross-validation with a maximum number of $15$ eigenfunctions and $9$ neighbors, respectively. In this case, the centroid method with the $FM_C$ semi-distance is the winner. Note that this method coincides in this case with the functional linear Bayes classification rule. The highest proportion of correct classification for the Phoneme dataset is $0.8238$ which is slightly larger than other alternatives.

On the other hand, Table $14$ shows the means and the standard deviation (between parentheses) of the number of principal components needed to calculate the $FPC_C$, $FPC_D$, $FM_C$, $FM_D$ and $DH$ semi-distances and the $FLBCR$ and $FQBCR$ methods for the Phoneme dataset. The mean numbers of functional principal components used with the winner methods is around $9$. However, other methods with worst performance have also mean values close to $9$. Therefore, in this case, the differences between performances are apparently due to the methods themselves.

\begin{table}
\begin{center}
\caption{Proportion of correct classification for the Phoneme dataset}
\begin{tabular}{cccccccccc}
\hline
& $L^{1}$ & $L^{2}$ & $L^{\infty }$ & $FPC_C$ & $FPC_{D}$ & $FM_C$ & $%
FM_{D}$ & $DH$ & $-$ \\ \hline
\textbf{kNN} & $\underset{\left( .0235\right) }{.7918}$ & $\underset{\left(
.0248\right) }{.7847}$ & $\underset{\left( .0258\right) }{.7838}$ & $%
\underset{\left( .0240\right) }{.7996}$ & $\underset{\left( .0233\right) }{%
.7799}$ & $\underset{\left( .0218\right) }{.8124}$ & $\underset{\left(
.0233\right) }{.7961}$ & $-$ & $-$ \\
\textbf{Centroid} & $\underset{\left( .0319\right) }{.7542}$ & $\underset{%
\left( .0307\right) }{.7386}$ & $\underset{\left( .0283\right) }{.7038}$ & $%
\underset{\left( .0307\right) }{.7401}$ & $\underset{\left( .0303\right) }{%
.7346}$ & $\underset{\left( .0236\right) }{\mathbf{.8238}}$ & $\underset{%
\left( .0218\right) }{.7994}$ & $\underset{\left( .0281\right) }{.8001}$ & $-
$ \\
\textbf{FLBCR} & $-$ & $-$ & $-$ & $-$ & $-$ & $-$ & $-$ & $-$ & $\underset{%
\left( .0236\right) }{\mathbf{.8238}}$ \\
\textbf{FQBCR} & $-$ & $-$ & $-$ & $-$ & $-$ & $-$ & $-$ & $-$ & $\underset{%
\left( .0218\right) }{.7994}$ \\
\textbf{LBCR Coef.} & $-$ & $-$ & $-$ & $-$ & $-$ & $-$ & $-$ & $-$ & $%
\underset{\left( .0250\right) }{.8050}$ \\
\textbf{QBCR Coef.} & $-$ & $-$ & $-$ & $-$ & $-$ & $-$ & $-$ & $-$ & $%
\underset{\left( .0261\right) }{.7802}$ \\ \hline
\end{tabular}
\end{center}
\end{table}

\begin{table}
\begin{center}
\caption{Means and standard deviations of the number of principal components used by the $FPC_C$, $FPC_D$, $FM_C$, $FM_D$ and $DH$ semi-distances and the $FLBCR$ and $FQBCR$ methods with the Phoneme dataset}
\begin{tabular}{ccccccc}
\hline
& $FPC_C$ & $FPC_{D}$ & $FM_C$ & $FM_{D}$ & $DH$ & $-$ \\ \hline
\textbf{kNN} & $\underset{\left( 2.87\right) }{8.31}$ & $\underset{\left(
3.26\right) }{9.37}$ & $\underset{\left( 2.82\right) }{9.08}$ & $\underset{%
\left( 3.40\right) }{9.48}$ & $-$ & $-$ \\
\textbf{Centroid} & $\underset{\left( 2.79\right) }{6.83}$ & $\underset{%
\left( 3.00\right) }{8.68}$ & $\underset{\left( 1.83\right) }{8.94}$ & $%
\underset{\left( 3.49\right) }{8.04}$ & $\underset{\left( 2.33\right) }{8.08}
$ & $-$ \\
\textbf{FLBCR} & $-$ & $-$ & $-$ & $-$ & $-$ & $\underset{\left( 1.83\right) }{%
8.94}$ \\
\textbf{FQBCR} & $-$ & $-$ & $-$ & $-$ & $-$ & $\underset{\left( 3.49\right) }{%
8.04}$ \\ \hline
\end{tabular}
\end{center}
\end{table}

\section{Conclusions}

This paper has introduced a new semi-distance for functional data that generalize the multivariate Mahalanobis distance to the functional framework. For that, it is used the regularized square root inverse operator given in Mas (2007) that allows to write the functional Mahalanobis semi-distance between an observation and the sample mean function of the set of functions in terms of the standardize functional principal component scores. Afterwards, new
versions of several classification procedures have been proposed based on the functional Mahalanobis semi-distance. Some Monte Carlo experiments and the analysis of two real data examples illustrate the good behavior of the
classification methods based on the functional Mahalanobis semi-distance. As mentioned previously, the range of applications of the functional Mahalanobis semi-distance is large and includes clustering, hypothesis testing and
outlier detection, among others. This would be the objective of future work.

\section*{Acknowledgements}

Financial support by MEC project ECO2012-38442 is gratefully acknowledged. The authors would like to thank Hugh Chipman for very helpful comments.

\section*{Appendix}

\subsection*{Proof of Proposition 2.1}

From (\ref{pro12}), it is possible to write:

\begin{gather*}
d_{FM}^{K}(\chi ,\mu _{\chi })=\left\langle \Gamma _{K}^{-\frac{1}{2}}(\chi -\mu _{\chi }),\Gamma _{K}^{-\frac{1}{2}}(\chi -\mu _{\chi
})\right\rangle ^{1/2}= \\
=\left\langle \sum_{k=1}^{K}\frac{1}{\lambda _{k}^{1/2}}(\psi _{k}\otimes\psi _{k})(\chi -\mu _{\chi }),\sum_{k=1}^{K}\frac{1}{\lambda _{k}^{1/2}}(\psi _{k}\otimes \psi _{k})(\chi -\mu _{\chi })\right\rangle ^{1/2}.
\end{gather*}
Now, from (\ref{prod1}) and (\ref{loeve23}), the previous expression leads to:
\begin{equation*}
d_{FM}^{K}(\chi ,\mu _{\chi })=\left\langle \sum_{k=1}^{K}\frac{1}{\lambda _{k}^{1/2}}\left[ \left\langle \psi _{k},\sum_{j=1}^{\infty }\theta
_{j}\psi _{j}\right\rangle \psi _{k}\right] ,\sum_{k=1}^{K}\frac{1}{\lambda_{k}^{1/2}}\left[ \left\langle \psi _{k},\sum_{j=1}^{\infty }\theta _{j}\psi
_{j}\right\rangle \psi _{k}\right] \right\rangle ^{1/2}.
\end{equation*}%
As the inner product is linear for $\theta _{k}$ and the $\psi _{k}$ are orthonormal eigenfunctions, it is possible to write:
\begin{equation*}
d_{FM}^{K}(\chi ,\mu _{\chi })=\left\langle \sum_{k=1}^{K}\frac{\theta
_{k}}{\lambda _{k}^{1/2}}\psi _{k},\sum_{k=1}^{K}\frac{\theta _{k}}{\lambda
_{k}^{1/2}}\psi _{k}\right\rangle ^{1/2}=\left(\sum_{k=1}^{K}\frac{s_{k}^{2}}{\lambda
_{k}}\right)^{1/2}=\left(\sum_{k=1}^{K}\omega_{k}^{2}\right)^{1/2}.
\end{equation*}

\subsection*{Proof of Proposition 2.2}

By hypothesis, the two functions $\chi _{1}$ and $\chi _{2}$ have the same mean function, $\mu _{\chi }$, and the same covariance operator, $\Gamma _{\chi }$. Therefore, from the Karhunen-Lo\`{e}ve expansion:
\begin{equation*}
\chi _{1}=\mu _{\chi }+\sum_{k=1}^{\infty }\theta _{1k}\psi _{k},
\end{equation*}%
and,
\begin{equation*}
\chi _{2}=\mu _{\chi }+\sum_{k=1}^{\infty }\theta _{2k}\psi _{k},
\end{equation*}
where $\theta _{1k}=\left\langle \chi _{1}-\mu _{\chi },\psi_{k}\right\rangle $ and $\theta _{2k}=\left\langle \chi _{2}-\mu _{\chi},\psi _{k}\right\rangle $, for $k=1,\ldots$ are the functional principal component scores of
$\chi _{1}$ and $\chi _{2},$ respectively. Consequently, the difference between the two functions $\chi _{1}$ and $\chi _{2}$ can be written as:
\begin{equation}
\chi _{1}-\chi _{2}=\sum_{k=1}^{\infty }(\theta _{1k}-\theta _{2k})\psi _{k}.\label{ael}
\end{equation}%
Using the expression (\ref{pro12}) of the regularized square root inverse operator, the Mahalanobis semi-distance between $\chi _{1}$ and $\chi _{2}$ is given by:
\begin{gather*}
d_{FM}^{K}(\chi _{1},\chi _{2})=\left\langle \Gamma _{K}^{-\frac{1}{2}}(\chi _{1}-\chi _{2}),\Gamma _{K}^{-\frac{1}{2}}(\chi _{1}-\chi_{2})\right\rangle ^{1/2}= \\
=\left\langle \sum_{k=1}^{K}\frac{1}{\lambda _{k}^{1/2}}(\psi _{k}\otimes\psi _{k})(\chi _{1}-\chi _{2}),\sum_{k=1}^{K}\frac{1}{\lambda _{k}^{1/2}}(\psi _{k}\otimes \psi _{k})(\chi _{1}-\chi _{2})\right\rangle ^{1/2}
\end{gather*}
Now, from (\ref{prod1}) and (\ref{ael}), the above expression can be written as:
\begin{gather*}
d_{FM}^{K}(\chi _{1},\chi _{2})=\left\langle \sum_{k=1}^{K}\frac{1}{\lambda _{k}^{1/2}}\left\langle \psi _{k},\chi _{1}-\chi _{2}\right\rangle\psi _{k},\sum_{k=1}^{K}\frac{1}{\lambda _{k}^{1/2}}\left\langle \psi
_{k},\chi _{1}-\chi _{2}\right\rangle \psi _{k}\right\rangle ^{1/2} = \\
=\sum_{k=1}^{K}\frac{1}{\lambda _{k}}\left\langle \left\langle \psi_{k},\sum_{j=1}^{\infty }(\theta _{1j}-\theta _{2j})\psi _{j}\right\rangle\psi _{k},\left\langle \psi _{k},\sum_{j=1}^{\infty }(\theta _{1j}-\theta
_{2j})\psi _{j}\right\rangle \psi _{k}\right\rangle ^{1/2} = \\
=\left(\sum_{k=1}^{K}\frac{1}{\lambda _{k}}\left\langle (\theta _{1k}-\theta_{2k})\psi _{k},(\theta _{1k}-\theta _{2k})\psi _{k}\right\rangle\right)^{1/2}=\left(\sum_{k=1}^{K}(\omega_{1k}-\omega_{2k})^{2}\right)^{1/2}
\end{gather*}
where $\omega_{1k}=\theta _{1k}/\lambda _{k}^{1/2}$ and $\omega_{2k}=\theta_{2k}/\lambda _{k}^{1/2}$, for $k=1,2,\ldots $ are the standardized functional principal component scores of $\chi _{1}$ and $\chi _{2}$,
respectively.

\subsection*{Proof of Proposition 2.3}

The proof of this proposition is trivial in view of Proposition 2.2 that asserts that $d_{FM}^{K}(\chi _{1},\chi _{2})$ is just the Euclidean distance between the first $K$ standardized functional principal component scores of
$\chi _{1}$ and $\chi _{2}$. Note that $d_{FM}^{K}(\chi _{1},\chi _{2})$ is not a functional distance because $d_{FM}^{K}(\chi _{1},\chi _{2})=0$ if $\chi _{1}$ and $\chi _{2}$ have the same first $K$ functional principal
component scores, which does not imply $\chi _{1}=\chi _{2}$.

\subsection*{Proof of Theorem 2.1}

The functional Mahalanobis semi-distance between the Gaussian process $\chi $ and its mean function $\mu _{\chi }$ is given in (\ref{dmk}). Now, as $\chi $ is a Gaussian process, the standardized functional principal component
scores, $\omega_{k}$, for $k=1,2,\ldots $ are independent standard Gaussian random variables (see, Ash and Gardner, 1975) that shows the result.

\section*{References}

Alonso, A. M., Casado, D. and Romo, J. (2012) Supervised classification for functional data: a weighted distance approach. Comput. Statist. Data Anal., 56, 2334-2346.

Araki, Y., Konishi, S., Kawano, S. and Matsui, H. (2009) Functional logistic discrimination via regularized basis expansions. Commun. Statist. Theory. Math., 38, 2944-2957.

Ash, R. B. and Gardner, M. F. (1975) \textit{Topics in stochastic processes}. Academic Press, New York.

Ba\'illo, A., Cuevas, A. and Cuesta-Albertos, J. A. (2011) Supervised classification for a family of Gaussian functional models. Scand. J. of Statist., 38, 480-498.

Biau, G., Bunea, F. and Wegkamp, M. H. (2005) Functional classification in Hilbert spaces. IEEE Trans. Inform. Theory., 51, 2163-2172.

C\'erou, F. and Guyader, A. (2006) Nearest neighbor classification in infinite dimension. ESAIM: Probability and Statistics, 10, 340-355.

Cuevas, A., Febrero, M. and Fraiman, R. (2007) Robust estimation and classification for functional data via projection-based depth notions. Comput. Statist., 22, 481, 496.

Delaigle, A. and Hall, P. (2012) Achieving Near Perfect Classification for Functional Data. J. R. Statist. Soc. B, 74, 267-286.

Epifanio, I. (2008) Shape descriptors for classification of functional data. Technometrics, 50, 284-294.

Ferraty, F. and Vieu, P. (2003) Curves Discrimination: A Nonparametric Functional Approach. Comput. Statist. Data Anal., 51, 4878-4890.

Ferraty, F. and Vieu, P. (2006) \textit{Nonparametric Functional Data Analysis}. Springer, New York.

Glendinning, R. H. and Herbert, R. A. (2003) Shape classification using smooth principal components. Pattern Recogn. Lett., 24, 2021-2030.

Hall, P., Poskitt, D. and Presnell, B. (2001) A functional data-analytic approach to signal discrimination. Technometrics, 43, 1-9.

Hall, P. and Hosseini-Nasab, M. (2006) On properties of functional principal components analysis. J. R. Statist. Soc. B, 68, 109-126.

James, G. M. and Hastie, T. J. (2001) Functional linear discriminant analysis for irregularly sampled curves. J. R. Statist. Soc. B, 63, 533-550.

Leng, X. Y. and M\"{u}ller, H.-G. (2006) Classification using functional data analysis for temporal gene expression data. Bioinformatics, 22, 68-76.

Li, B. and Yu, Q. (2008) Classification of Functional Data: A segmentation Approach. Comput. Statist. Data Anal., 52, 4790-4800.

L\'opez-Pintado, S. and Romo, J. (2006) \textit{Depth-based classification for functional data}. In DIMACS Series in Discrete Mathematics and Theoretical Computer Science, vol. 72, 103-120. Providence: American Mathematical Society.

Mahalanobis, P. C. (1936) On the generalized distance in Statistics. Proc. Natl. Acad. Sci. India, 12, 49-55.

Martin-Barragan, B., Lillo, R. E. and Romo, J. (2013) Interpretable support vector machines for functional data. European J. of Oper. Res., In press.

Mas, A. (2007) Weak convergence in the functional autoregressive model. J. Mult. Anal., 98, 1231-1261.

Preda, C., Saporta, G. and Leveder, C. (2007) PLS classification of functional data. Comput. Statist., 22, 223-235.

Ramsay, J. O. and Silverman, B. W. (2005) \textit{Functional Data Analysis}. 2nd Edition. Springer, New York.

Rossi, F. and Villa N. (2006) Support Vector Machine for Functional Data Classification. Neurocomputing, 69, 730-742.

Shin, H. (2008) An extension of Fisher´s discriminant analysis for stochastic processes. J. Multiv. Anal., 99, 1191-1216.

Song, J. J., Deng, W., Lee, H.-J. and Kwon, D. (2008) Optimal classification for time-course gene expression data using functional data analysis. Comput Biol. Chem., 32, 426-432.

Wang, X. H., Ray, S. and Mallick, B. K. (2007) Bayesian curve classification using wavelets. J. Am. Statist. Ass., 102, 962-973.

\end{document}